\documentclass[11pt]{article}
\usepackage{t1enc}
\usepackage{lmodern}
\usepackage[T2A,T1]{fontenc}
\usepackage[utf8]{inputenc}
\usepackage[russian,english]{babel}
\usepackage{amsmath,amssymb,amsfonts,amsthm,mathrsfs,textcomp,url,bbm,cancel,comment}
\usepackage{graphicx}
\usepackage{float}
\usepackage[all,ps]{xy}
\usepackage[colorlinks,urlcolor=cyan,citecolor=blue,linkcolor=blue]{hyperref}
\usepackage[dvipsnames]{xcolor}
\usepackage[margin=3.5cm]{geometry}
\usepackage{longtable}
\usepackage{tikz}
\DeclareRobustCommand{\cyrins}[1]{%
  \begingroup\fontfamily{cmr}%
  \foreignlanguage{russian}{#1}%
  \endgroup}
\usetikzlibrary{matrix,calc,decorations,positioning}
\pgfkeys{/tikz/.cd,
    alt double distance/.initial=5pt,
    alt double step/.initial=1pt,}
\pgfdeclaredecoration{double deco}{initial}
{
\state{initial}[width=\pgfkeysvalueof{/tikz/alt double step},next state=cont] {
    \pgfmoveto{\pgfpoint{\pgfkeysvalueof{/tikz/alt double step}}{\pgfkeysvalueof{/tikz/alt double distance}/2}}
    \pgfpathlineto{\pgfpoint{0.3\pgflinewidth}{\pgfkeysvalueof{/tikz/alt double distance}/2}}
    \pgfpathmoveto{\pgfpoint{0.3\pgflinewidth}{-\pgfkeysvalueof{/tikz/alt double distance}/2}}
    \pgfpathlineto{\pgfpoint{1pt}{-\pgfkeysvalueof{/tikz/alt double distance}/2}}
    \pgfcoordinate{lastup}{\pgfpoint{1pt}{\pgfkeysvalueof{/tikz/alt double distance}/2}}
    \pgfcoordinate{lastdown}{\pgfpoint{1pt}{-\pgfkeysvalueof{/tikz/alt double distance}/2}}}
  \state{cont}[width=\pgfkeysvalueof{/tikz/alt double step}]{
     \pgfmoveto{\pgfpointanchor{lastup}{center}}
     \pgfpathlineto{\pgfpoint{\pgfkeysvalueof{/tikz/alt double step}}{\pgfkeysvalueof{/tikz/alt double distance}/2}}
     \pgfcoordinate{lastup}{\pgfpoint{\pgfkeysvalueof{/tikz/alt double step}}{\pgfkeysvalueof{/tikz/alt double distance}/2}}
     \pgfmoveto{\pgfpointanchor{lastdown}{center}}
     \pgfpathlineto{\pgfpoint{\pgfkeysvalueof{/tikz/alt double step}}{-\pgfkeysvalueof{/tikz/alt double distance}/2}}
     \pgfcoordinate{lastdown}{\pgfpoint{\pgfkeysvalueof{/tikz/alt double step}}{-\pgfkeysvalueof{/tikz/alt double distance}/2}}}
  \state{final}[width=0pt]
  { 
    \pgfmoveto{\pgfpointdecoratedpathlast}}}
\tikzset{alt double/.style={decorate,decoration=double deco}}
\allowdisplaybreaks

\newcommand{\toverset}[2]{%
\mathop{#2}\limits^{\vbox to -.1ex{\kern-0.4ex\hbox{$\scriptstyle #1$}\vss}}}
\newcommand{\tightoverset}[2]{%
\mathop{#2}\limits^{\vbox to -.5ex{\kern-0.4ex\hbox{$\scriptstyle #1$}\vss}}}
\renewcommand{\underset}[2]{%
\mathop{#2}\limits_{\vbox to -.5ex{\kern-1.6ex\hbox{$\scriptstyle #1$}\vss}}}
\newcommand{\tightunderset}[2]{%
\mathop{#2}\limits_{\vbox to -.5ex{\kern-1.8ex\hbox{$\scriptstyle #1$}\vss}}}
\newcommand{\xra}[1]{%
\mathop{\xrightarrow{~#1~}}}

\newcommand{\cpt}[1]{%
\mathop{\toverset{^\bullet}{#1}}}

\newcommand{\usqcup}[1]{%
\mathop{\tightunderset{#1}{\sqcup}}}
\newcommand{\utimes}[1]{%
\mathop{\tightunderset{#1}{\times}}}

\def\N{\mathbb{N}}
\def\Z{\mathbb{Z}}
\def\Q{\mathbb{Q}}
\def\R{\mathbb{R}}

\def\RP{\mathbb{R}P}
\def\CP{\mathbb{C}P}

\def\AA{\mathscr{A}}
\def\CC{\mathscr{C}}
\def\FF{\mathscr{F}}
\def\GG{\mathscr{G}}
\def\VV{\mathscr{V}}
\def\NN{\mathfrak{N}}
\def\WW{\mathfrak{W}}
\def\onto{\twoheadrightarrow}
\def\into{\hookrightarrow}
\def\imto{\looparrowright}

\def\ol{\overline}

\def\imm{\mathrm{imm}}
\def\fold{\mathrm{fold}}
\def\im{\mathop{\rm im}}
\def\coker{\mathop{\rm coker}}

\def\codim{\mathop{\rm codim}}

\def\id{\mathop{\rm id}}
\def\pr{\mathop{\rm pr}}

\def\PD{\mathop{\rm PD}}

\def\O{\mathrm{O}}
\def\SO{\mathrm{SO}}

\def\rk{\mathop{\rm rk}}
\def\i{\mathrm{int}}
\def\fr{\mathrm{fr}}

\def\Emb{\textstyle{\mathop{\rm Emb}}}

\def\Imm{\textstyle{\mathop{\rm Imm}}}
\def\Fold{\textstyle{\mathop{\rm Fold}}}
\def\Mor{\textstyle{\mathop{\rm Mor}}}
\def\Cob{\textstyle{\mathop{\rm Cob}}}

\renewcommand{\Lambda}{\det}
\newenvironment{prf}%
{\par\noindent\textbf{Proof.\enspace\ignorespaces}}%
{~$\square$\par}%
\newenvironment{sclaim}%
{\par\noindent\textit{Claim.\enspace\ignorespaces}\begin{em}}%
{\end{em}\par}%
\newenvironment{sprf}%
{\medskip\par\noindent\textit{Proof.\enspace\ignorespaces}}%
{~$\diamond$\par\medskip}%
{\medskip\par\noindent\textit{Remark.\enspace\ignorespaces}}%
{\par\medskip}%
\newenvironment{key}%
{\par\noindent\begin{small}\textbf{Keywords.\enspace\ignorespaces}}%
{\end{small}\par}%
\newenvironment{ack}%
{\par\noindent\begin{small}\textbf{Acknowledgement.\enspace\ignorespaces}}%
{\end{small}\par}%
\newtheorem{thm}{Theorem}[section]%
\newtheorem{lemma}[thm]{Lemma}%
\newtheorem{prop}[thm]{Proposition}%
\newtheorem{crly}[thm]{Corollary}%
\theoremstyle{definition}
\newtheorem{defi}[thm]{Definition}%
\newtheorem{rmk}[thm]{Remark}%
\newtheorem{ex}[thm]{Example}%
\newcommand{\ZeroRoman}[1]{%
\ifcase\value{#1}\relax 0\else\Roman{#1}\fi}

\newcounter{t}
\renewcommand\thet{\Roman{t}}
\numberwithin{equation}{section}
\title{Analogues of the Atiyah--Wall exact sequences\\for cobordism groups of singular maps}
\author{András Csépai\thanks{ELTE Eötvös Loránd University, Budapest, Hungary}}
\date{}
\begin{document}

\maketitle

\begin{abstract}
Classical results of Rohlin, Dold, Wall and Atiyah yield two exact sequences that connect the oriented and unoriented (abstract) cobordism groups $\Omega_n$ and $\NN_n$. In this paper we present analogous exact sequences connecting the oriented and unoriented cobordism groups of maps with prescribed singularities. This gives positive answer to a fifteen-year-old question posed by Szűcs and has interesting consequences even in the case of cobordisms of immersions.
\end{abstract}

\medskip

\begin{key}
cobordism; singular maps; exact sequences
\end{key}

\setcounter{part}{-1}
\part{Introduction}

\section{Classical results}

In the 1950's much work was done on the determination of the structures and generators of the oriented and unoriented abstract cobordism rings $\Omega_*$ and $\NN_*$. One of the main tools for these computations was the existence of two exact sequences which resulted from individual works of Rohlin \cite{rohlin}, Dold \cite{dold}, Wall \cite{wallcob} and finally from a conceptual method by Atiyah \cite{atiyah}.

Recall that the elements of $\Omega_n$ are represented by oriented closed smooth ($C^\infty$) $n$-manifolds and two such manifolds are cobordant if they together bound an oriented compact $(n+1)$-manifold with boundary with matching orientations; the elements of $\NN_n$ are similar but without orientations. Between these two groups stands the cobordism group $\WW_n$ of the so-called \textit{Wall manifolds} (see e.g. \cite{stong}) i.e. closed manifolds whose first Stiefel--Whitney class is the mod $2$ reduction of an integer cohomology class; two such manifolds are cobordant if they together bound a compact $(n+1)$-manifold with boundary with matching integer first Stiefel--Whitney classes. Now the forgetful homomorphism from $\Omega_n$ to $\NN_n$ (i.e. the map we get by ignoring the orientation) is the composition $\Omega_n\to\WW_n\to\NN_n$ of two forgetful homomorphisms.

The classical Atiyah--Wall exact sequences (see \cite[theorems 4.2 and 4.3]{atiyah}) are a long exact sequence
\begin{equation}\label{ces1}
\ldots\to\Omega_n\to\Omega_n\to\WW_n\to\Omega_{n-1}\to\ldots
\tag{\text{I}}
\end{equation}
and a short exact sequence
\begin{equation}\label{ces2}
0\to\WW_n\to\NN_n\to\NN_{n-2}\to0
\tag{\text{II}}
\end{equation}
containing these forgetful homomorphisms.

The main result of the present paper is the generalisation of these sequences to the cobordism theory of singular maps which answers an open question of Szűcs proposed in \cite[section 19]{hosszu}. The precise statements of this result can be found in theorems \ref{t1} and \ref{t2} at the end of the next section after recalling and introducing the necessary definitions. Theorem \ref{t1} and theorem \ref{t2} generalise the sequences (\ref{ces1}) and (\ref{ces2}) respectively to cobordism groups of singular maps and they will be proved in part \ref{p1} and part \ref{p2} respectively; then in part \ref{p3} we shall apply them to compute various cobordism groups and finally to obtain analogous exact sequences for bordism groups.

\section{Cobordism groups of singular maps}

Throughout this paper we consider smooth ($C^\infty$) maps of $n$-manifolds to $(n+k)$-manifolds where $k$ is a fixed positive integer and $n$ is arbitrary. If we want to indicate the dimension of a manifold, we put it in a superindex (i.e. $M^n$ means that $M$ is a manifold of dimension $n$) but in most cases we omit this index. If we do not state otherwise, then we also always make the technical assumption that any smooth map between manifolds $f\colon M\to P$ is such that $f^{-1}(\partial P)=\partial M$ and $f$ is transverse to the boundary $\partial P$. 

\begin{defi}
Two smooth map germs
$$\eta\colon(\R^n,0)\to(\R^{n+k},0)\quad\text{and}\quad\vartheta\colon(\R^n,0)\to(\R^{n+k},0)$$
are said to be \textit{left-right equivalent} (in some sources also called \textit{$\AA$-equivalent}) if there are diffeomorphism germs $\varphi$ and $\psi$ of $(\R^n,0)$ and $(\R^{n+k},0)$ respectively such that $\vartheta=\psi\circ\eta\circ\varphi^{-1}$. The \textit{suspension} of the germ $\eta$ is the germ
$$\eta\times{\id}_{\R^1}\colon(\R^{n+1},0)\to(\R^{n+k+1},0).$$
By the {\it singularity class} (or simply {\it singularity}) of $\eta$ we mean the equivalence class of $\eta$ in the equivalence relation generated by left-right equivalence and suspension. The singularity class of $\eta$ is denoted by $[\eta]$.
\end{defi}

Observe that in each singularity class the codimension $k$ of the germs is fixed but the dimension $n$ is not.

\begin{defi}
For manifolds $M^n,P^{n+k}$ the product of the diffeomorphism groups of $M$ and $P$ has an action on the space $C^\infty(M,P)$ of smooth maps $f\colon M\to P$ defined by the formula $(\varphi,\psi)\mapsto\psi\circ f\circ\varphi^{-1}$. Two maps $f,g\colon M\to P$ are said to be \textit{left-right equivalent} if they are in the same orbit of this action, i.e. $g=\psi\circ f\circ\varphi^{-1}$ for some diffeomorphisms $\varphi$ of $M$ and $\psi$ of $P$. A map $f\colon M\to P$ is said to be \textit{stable} if it is in the interior of an orbit, i.e. it has a neighbourhood $U\subset C^\infty(M,P)$ such that every element of $U$ is left-right equivalent to $f$. We define a map germ a \textit{stable germ} if it is the germ of a stable map.
\end{defi}

In the present paper we only consider stable germs. Note that if the germ $\eta$ is stable, then every germ representing the singularity $[\eta]$ is also stable, hence it is justified to call singularities of stable germs \textit{stable singularities}.

\begin{rmk}
The restriction to only study maps with stable germs is quite mild, in the so-called nice dimensions stable maps form an open dense subset in the space of smooth maps; see \cite{math}.
\end{rmk}

Let us now fix a set $\tau$ of singularities of stable $k$-codimensional germs.

\begin{defi}
A smooth map $f\colon M^n\to P^{n+k}$ is said to be a {\it $\tau$-map} if all of its germs belong to singularity classes in $\tau$. For a singularity class $[\eta]\in\tau$ a point $p\in M$ is said to be an {\it $[\eta]$-point} if the germ of $f$ at $p$ is equivalent to $\eta$; the set of $[\eta]$-points in $M$ is denoted by $[\eta](f)$. By a slight abuse of notation we will sometimes write $\eta(f)$ instead of $[\eta](f)$ for better readability.
\end{defi}

\begin{defi}
For a singularity class $[\eta]$ the minimal number $m$ for which there is a germ $\vartheta\colon(\R^m,0)\to(\R^{m+k},0)$ in the singularity class $[\eta]$ is called the {\it codimension of the singularity} $[\eta]$ and is denoted by $\codim[\eta]$.
\end{defi}

\begin{rmk}
If $f\colon M\to P$ is a $\tau$-map, then $M$ is naturally stratified by the submanifolds $\eta(f)$ for all $[\eta]\in\tau$. Here $\eta(f)$ is a submanifold of $M$ of codimension $\codim[\eta]$.
\end{rmk}

\begin{ex}
~
\begin{enumerate}
\item If $\tau$ only contains the class $\Sigma^0$ of regular germs (i.e. those with maximal-rank derivative at $0$), then $\tau$-maps are just the immersions.
\item If $\tau$ consists of $\Sigma^0$ and all stable singularities of type $\Sigma^1$, that is, singularities with representatives whose derivative at $0$ has corank $1$, then $\tau$-maps are called \textit{Morin maps}. Further restricting the set $\tau$ of allowed singularities we can obtain e.g. the so-called fold maps (i.e. $\{\Sigma^0,\Sigma^{1,0}\}$-maps) and cusp maps (i.e. $\{\Sigma^0,\Sigma^{1,0},\Sigma^{1,1,0}\}$-maps) where fold ($\Sigma^{1,0}$) and cusp ($\Sigma^{1,1,0}$) are the two simplest types of Morin singularities; see \cite{mor}.
\end{enumerate}
\end{ex}

In the following we shall work with $\tau$-maps to a fixed target manifold $P$ up to a cobordism relation defined as follows:

\begin{defi}\label{cob}
We call two $\tau$-maps $f_0\colon M_0^n\to P^{n+k}$ and $f_1\colon M_1^n\to P^{n+k}$ (with closed source manifolds $M_0$ and $M_1$) {\it $\tau$-cobordant} if there is
\begin{itemize}
\item[(i)] a compact manifold $W^{n+1}$ with boundary such that $\partial W=M_0\sqcup M_1$,
\item[(ii)] a $\tau$-map $F\colon W^{n+1}\to P\times[0,1]$ such that for $i=0,1$ we have $F^{-1}(P\times\{i\})=M_i$ and $F|_{M_i}=f_i$.
\end{itemize}
The $\tau$-cobordism class of $f\colon M^n\to P^{n+k}$ is denoted by $[f]$ and the set of all $\tau$-cobordism classes of $\tau$-maps to the manifold $P$ is denoted by $\Cob_\tau(P)$.
\end{defi}

The set $\Cob_\tau(P)$ admits a natural commutative semigroup operation by the disjoint union: if $f\colon M^n\to P^{n+k}$ and $g\colon N^n\to P^{n+k}$ are $\tau$-maps, then so is
$$f\sqcup g\colon M\sqcup N\to P$$
and the $\tau$-cobordism class $[f]+[g]:=[f\sqcup g]$ is well-defined. This operation has a null element represented by the empty map, moreover, it is actually a group operation (the inverse of any element in $\Cob_\tau(P)$ is explicitly constructed in \cite{hosszu}). This way $\Cob_\tau(P)$ becomes an Abelian group for any manifold $P$.

Next we shall endow $\tau$-maps with various stable normal structures. These will be defined in the following four points; loosely speaking they are decorations on the stable normal bundle, that is, for a map $f\colon M\to P$ the stable isomorphism type of the virtual vector bundle $\nu_f:=f^*TP\ominus TM$. If $\sigma$ is such a decoration, then $\tau$-maps equipped with $\sigma$-structures will be denoted $\tau^\sigma$-maps. In this paper we shall work with the following types of $\tau^\sigma$-maps:
\begin{enumerate}
\item $\sigma=G$: Let $G$ be a stable group using the definition of Wall \cite[section 8.2]{stabgr} which we recall now:

  \begin{defi}
    \label{stabgr}
    A {\it stable group} $G$ is defined as the direct limit of a sequence of group homomorphisms $\iota_n\colon G(n)\to G(n+1)$ where
    \begin{itemize}
    \item[(i)] there are homomorphisms $\alpha_n\colon G(n)\to\O(n)$ such that the diagram
      $$\xymatrix{
      \ldots\ar[r]^(.45){\iota_{n-1}} & G(n)\ar[d]_{\alpha_n}\ar[r]^(.4){\iota_n} & G(n+1)\ar[d]^{\alpha_{n+1}}\ar[r]^(.6){\iota_{n+1}} & \ldots\\
      \ldots\ar@{^(->}[r] & \O(n)\ar@{^(->}[r] & \O(n+1)\ar@{^(->}[r] & \ldots
      }$$
      is commutative where the inclusions in the lower row are natural,
    \item[(ii)] there is a weakly increasing function $c\colon\N\to\N$ tending to $\infty$ such that the map $\iota_n$ is $c(n)$-connected,
    \item[(iii)] there are homomorphisms $\beta_{n,m}\colon G(n)\times G(m)\to G(n+m)$ such that the diagrams
      $$\xymatrix@C=3pc{
        G(n)\times G(m)\ar[d]_{\iota_n\times\id}\ar[r]^{\beta_{n,m}} & G(n+m)\ar[d]^{\iota_{n+m}} & G(n)\times G(m)\ar[l]_{\beta_{n,m}}\ar[d]^{\id\times\iota_m} \\
        G(n+1)\times G(m)\ar[r]^{\beta_{n+1,m}} & G(n+m+1) & G(n)\times G(m+1)\ar[l]_{\beta_{n,m+1}}
      }$$
      $$\xymatrix{
        G(n)\times G(m)\ar[d]_{\alpha_n\times\alpha_m}\ar[r]^{\beta_{n,m}} & G(n+m)\ar[d]^{\alpha_{n+m}} \\
        \O(n)\times\O(m)\ar@{^(->}[r] & \O(n+m)
      }$$
      and
      $$\xymatrix@C=3.5pc{
        G(n)\times G(m)\times G(l)\ar[d]_{\beta_{n,m}\times\id}\ar[r]^{\id\times\beta_{m,l}} & G(n)\times G(m+l)\ar[d]^{\beta_{n,m+l}} \\
        G(n+m)\times G(l)\ar[r]^{\beta_{n+m,l}} & G(n+m+l)
      }$$
      commute up to conjugation by an element in the component of the identity,
    \item[(iv)] there is a commutative diagram
      $$\xymatrix{
        G(n)\times G(m)\ar[d]_{\iota_{n+m}\circ\beta_{n,m}}\ar[r]^{\gamma_{n,m}} & G(m)\times G(n)\ar[d]^{\iota_{n+m}\circ\beta_{m,n}} \\
        \O(n+m)\ar[r]^{\delta_{n,m}} & \O(n+m)
      }$$
      where $\gamma_{n,m}$ denotes the interchange of factors and $\delta_{n,m}$ is the conjugation by an element whose determinant has sign $(-1)^{nm}$.
    \end{itemize}
  \end{defi}

  We define {\it $\tau^G$-maps} to be $\tau$-maps $f\colon M^n\to P^{n+k}$ for which the virtual normal bundle $\nu_f$ has structure group $G$ in the following sense:

  \begin{defi}
    For a stable group $G$ and a virtual vector bundle $\nu$ we say that {\it $\nu$ has structure group $G$} if there is a fixed equivalence class of vector bundles stably isomorphic to $\nu$ such that the structure group of each vector bundle $\xi$ is reduced to $G(\rk\xi)$ and their equivalence is defined as follows. Note that the reduction of the structure group of $\xi$ to $G(\rk\xi)$ also gives the reduction of the structure group of $\xi\oplus\varepsilon^r$ to $G(\rk\xi+r)$ for any $r$. Here and later on $\varepsilon^r$ means the trivial rank-$r$ vector bundle over any base space. Now two such bundles $\xi,\zeta$ represetning $\nu$ are said to be equivalent if for some $r$ the bundles $\xi\oplus\varepsilon^{r-\rk\xi}$ and $\zeta\oplus\varepsilon^{r-\rk\zeta}$ are isomorphic as $G(r)$-bundles.
  \end{defi}

For example if $G=\SO$, then $\tau^\SO$-maps are $\tau$-maps with oriented normal bundles; if $G=\O$, then $\tau^\O$-maps are just $\tau$-maps without further conditions. If $\tau$ consists of all possible singularities of $k$-codimensional germs,
 then any manifold $M^n$ has a $\tau$-map to $\R^{n+k}$ uniquely up to $\tau$-cobordism and its stable normal bundle is just the inverse of the tangent bundle $TM$ in any $K$-group of $M$; the same is true for any abstract cobordism of manifolds, hence in this case we have
$$\Cob_\tau^\SO(\R^{n+k})=\Omega_n\quad\text{and}\quad\Cob_\tau^\O(\R^{n+k})=\NN_n.$$
As another example, if $\tau=\{\Sigma^0\}$ (i.e. $\tau$-maps are the immersions), then $\Cob_\tau^G(\R^{n+k})$ is the cobordism group of immersions of $n$-manifolds into $\R^{n+k}$ with normal structure group $G$ which was first described by Wells \cite{wells}. For various other singularity sets $\tau$ with $G$ being mostly $\SO$ or $\O$ the groups $\Cob_\tau^G(P)$ were considered e.g. by Szűcs \cite{analog}, \cite{hosszu}\footnote{We shall refer many times to \cite{hosszu} which only considers $\tau^\SO$-maps but contains theorems which work with the same proofs for more general stable normal structures as well, hence in this paper we will refer to them in their general form.}, Rimányi \cite{rsz}, Ando \cite{ando}, Terpai \cite{2k+2}, Kalmár \cite{kal} and Sadykov \cite{sad} (although in some of the papers cited the codimension $k$ of the maps is non-positive, unlike in the present paper).
\item $\sigma=\i$: We define {\it $\tau^\i$-maps} to be $\tau$-maps with the first Stiefel--Whitney class being an integer class, that is, a $\tau^\i$-map is a pair $(f,w)$ where $f\colon M^n\to P^{n+k}$ is a $\tau$-map and $w\in H^1(M;\Z)$ is a cohomology class such that its mod $2$ reduction is $w_1(\nu_f)$. Following e.g. Stong \cite{stong} we call this type of normal structure a {\it Wall structure} and also call $\tau^\i$-maps {\it Wall $\tau$-maps}. Again, if $\tau$ consists of all possible singularities of $k$-codimensional germs, then any manifold $M^n$ has a $\tau$-map to $\R^{n+k}$ uniquely up to $\tau$-cobordism and its first normal Stiefel--Whitney class is the same as $w_1(M)$, hence for this $\tau$ we have
$$\Cob_\tau^\i(\R^{n+k})=\WW_n.$$
\item $\sigma=\oplus m$ for a natural number $m$: {\it $\tau^{\oplus m}$-maps} are also called {\it $m$-framed $\tau$-maps} and they are defined so that their germs are of the form $\eta\times\id_{\R^m}$ (for $[\eta]\in\tau$) and the change of coordinate neighbourhoods always induces the identity on $\R^m$ here. This extends the notion of $m$-framed immersions, i.e. immersions equipped with $m$ pointwise independent normal vector fields; if $\tau=\{\Sigma^0\}$, then $\tau^{\oplus m}$-maps are just immersions with normal $m$-framing. We can endow $\tau^{\oplus m}$-maps with other stable normal structures $\sigma$ as well, and this defines {\it $\tau^{\sigma\oplus m}$-maps}. These maps were constructed in \cite[definition 9 and remark 10]{hosszu} and a very important property of them is the following: if the target manifold $P$ is of the form $Q\times\R^m$ for a manifold $Q$, then we have
\begin{equation}\label{+m}
\Cob_\tau^{\sigma\oplus m}(Q\times\R^m)\cong\Cob_\tau^\sigma(Q)
\end{equation}
and the isomorphism is given by a natural correspondence between framed and non-framed maps; see \cite[proposition 13]{hosszu}.
\item $\sigma=G\oplus\xi$: Let $G$ be a stable group and $\xi$ a vector bundle over the classifying space $BG$. Then {\it $\tau^{G\oplus\xi}$-maps} are generalisations of $\tau^{G\oplus m}$-maps and we postpone their precise definition to section \ref{tauxi}. Intuitively a $\tau^{G\oplus\xi}$-map is a $\tau^G$-map such that a bundle induced from $\xi$ splits off from its ``normal bundle'' (the sense in which we mean this will be clarified in section \ref{tauxi}); if $\xi$ is trivial of rank $m$, this just gives an $m$-framed $\tau^G$-map.
\end{enumerate}

Now if $\sigma$ is any of the above four stable normal structures, then the cobordism of two $\tau^\sigma$-maps to a manifold $P$ can be defined by adding the $\sigma$-structure to the map $F$ in definition \ref{cob} and the cobordism group of $\tau^\sigma$-maps to $P$ is denoted by $\Cob_\tau^\sigma(P)$.


To form the statements of our main theorems it remains to define a vector bundle which will play the role of $\xi$ in a special type of $\tau^{G\oplus\xi}$-maps. We put $G:=\O$ and denote by $\Lambda\gamma$ the line bundle over $B\O$ induced by the first Stiefel--Whitney class $w_1\colon B\O\to\RP^\infty$ from the tautological line bundle. The notation of this line bundle is due to the fact that if $\gamma_n^\O$ denotes the tautological rank-$n$ vector bundle over $B\O(n)$, then the restriction of $\Lambda\gamma$ over $B\O(n)$ is the determinant bundle $\Lambda\gamma_n^\O$. In our second theorem we shall use the rank-$2$ bundle $2\Lambda\gamma=\Lambda\gamma\oplus\Lambda\gamma$.

Now we are in a position to state our main results which are as follows:

\medskip
\noindent\textbf{Theorem \refstepcounter{t}\thet\label{t1}.\enspace\ignorespaces}\textit{
For any set $\tau$ of stable singularities and any manifold $Q^q$ there is a long exact sequence
\begin{alignat*}2
\ldots&\xra{\psi_{m+1}}\Cob_\tau^\SO(Q\times\R^m)\xra{\varphi_m}\Cob^\SO_\tau(Q\times\R^m)\xra{\chi_m}\Cob_\tau^\i(Q\times\R^m)\xra{\psi_m} \\
&\xra{\psi_m}\Cob_\tau^\SO(Q\times\R^{m-1})\to\ldots
\end{alignat*}
where $\varphi_m$ is of the form $\id+\iota$ where $\iota$ is the involution represented by the reflection to a hyperplane, $\chi_m$ is the forgetful homomorphism and $\psi_m$ assigns to a cobordism class $[(f,w)]$ the class $[f|_{\PD(w)}]$.
}\medskip

\noindent\textbf{Theorem \refstepcounter{t}\thet\label{t2}.\enspace\ignorespaces}\textit{
For any set $\tau$ of stable singularities and any manifold $Q^q$ there is a long exact sequence
\begin{alignat*}2
\ldots&\xra{\psi'_{m+1}}\Cob_\tau^\i(Q\times\R^m)\xra{\varphi'_m}\Cob^\O_\tau(Q\times\R^m)\xra{\chi'_m}\Cob_\tau^{\O\oplus2\Lambda\gamma}(Q\times\R^m)\xra{\psi'_m} \\
&\xra{\psi'_m}\Cob_\tau^\i(Q\times\R^{m-1})\to\ldots
\end{alignat*}
where $\varphi'_m$ is the forgetful homomorphism and $\chi'_m$ assigns to a cobordism class $[f]$ the class $[f|_{\PD(w_1(\nu_f)^2)}]$.
}

\begin{rmk}\label{psin0}
We shall give a geometric description of $\psi'_m$ later in remark \ref{psi'}. Here we just note that its classical analogue is always zero while in our general case this does not always hold; proposition \ref{n0} will show that even the composition $\psi_{m-1}\circ\psi'_m$ is not necessarily zero.
\end{rmk}

\begin{rmk}
Although we made assumptions on the set $\tau$ of allowed singularities which generally exclude the set of all singularities, i.e. bordism groups in general cannot be considered as $\tau$-cobordism groups, the above theorems still hold for bordism groups as well; see theorem \ref{b}.
\end{rmk}

\begin{rmk}
\label{cohrmk}
Recall that the (say oriented) abstract cobordism and bordism groups give rise to an extraordinary cohomology theory $Q\mapsto M\SO^*(Q)$; see \cite{atiyah}. Similarly cobordism groups of singular maps also yield cohomology theories: for the singularity sets $\tau$ and stable normal structures $\sigma$ used in the present paper we can define extraordinary cohomology functors $h^*_{\tau^\sigma}$ (see \cite[section 19]{hosszu}) and if $Q$ is a manifold, then we have
$$h^{-m}_{\tau^\sigma}(Q)=\Cob_\tau^\sigma(Q\times\R^m)\quad\text{and}\quad h^m_{\tau^\sigma}(Q)=\Cob_\tau^{\sigma\oplus m}(Q)$$
for any natural number $m$ (cf. the isomorphism (\ref{+m})). 
We will not use this fact in our proofs, however we note that the exact sequences claimed above are sequences of these cohomology groups, hence they naturally extend infinitely to the right if we put $\Cob_\tau^{\sigma\oplus-m}(Q)$ in the place of $\Cob_\tau^\sigma(Q\times\R^m)$ for negative $m$'s. 
\end{rmk}

Theorems \ref{t1} and \ref{t2} will be proved in part \ref{p1} and part \ref{p2} respectively and in both cases we begin by constructing classifying spaces necessary for the proofs, then we define the above long exact sequence of cobordism groups, and finish by showing that the homomophisms in the sequence are the ones we claimed. We shall also see that these exact sequences generalise the classical exact sequences (\ref{ces1}) and (\ref{ces2}); see remarks \ref{gen1} and \ref{gen2}. After this, in part \ref{p3} we will see applications of theorems \ref{t1} and \ref{t2} to cobordism groups of immersions and Morin maps, and finally, bordism groups.

Note that in the sequence (\ref{ces1}) the arrow $\Omega_n\to\Omega_n$ is the multiplication by $2$, hence for theorem \ref{t1} to be completely analogous to it we should have $\varphi_m=2\id$, i.e. $\iota=\id$, which is indeed so in some special cases but not generally as we shall see. Moreover, $\psi_m$ is $[(f,w)]\mapsto[f|_{\PD(w)}]$ and not $[(f,w)]\mapsto[f|_{\PD(w_1(\nu_f))}]$ although its classical version $\WW_n\to\Omega_{n-1}$ in (\ref{ces1}) is just $[M]\mapsto[\PD(w_1(M))]$. This is because $\WW_n$ embeds into $\NN_n$, hence taking the Poincaré dual of the integer first Stiefel--Whitney class of a Wall manifold $M$ up to cobordism is independent of which integer representative of $w_1(M)$ we take, however, remark \ref{psin0} implies that $\varphi'_m$ is not always mono and so this simplification which works for abstract cobordism groups may not work in the general case.

Theorem \ref{t2} is not a perfect analogue of the sequence (\ref{ces2}) either: it is actually analogous to \cite[theorem 4.3]{atiyah} (if we do not take into account the normal $\O\oplus2\Lambda\gamma$-structure which is meaningless in the case of abstract cobordism groups). Atiyah in his paper then proves that the long exact sequence splits to short exact sequences yielding (\ref{ces2}) but, as noted in remark \ref{psin0} above, this is not true generally.


\part{The oriented case}\label{p1}

\section{Cobordism of $\tau^\i$-maps}\label{tauw}

One of the most important ingredients of any type of cobordism theory is an analogue of the Pontryagin--Thom construction which assigns a classifying space to that theory and gives a bijection between the cobordism classes and the homotopy classes of maps to this space. The classifying space of the cobordism groups $\Cob_\tau^\sigma(P)$ (where $\tau$ is a set of stable singularities, $\sigma$ is a stable normal structure and $P$ is an arbitrary manifold) will be denoted by $X_\tau^\sigma$, that is, $X_\tau^\sigma$ is the (homotopically unique) space which has the property
$$\Cob_\tau^\sigma(P)\cong[\cpt P,X_\tau^\sigma]$$
where $\cpt P$ denotes the one-point compactification of $P$ and $[\cdot,\cdot]$ means the (based) homotopy classes of maps from the first space to the second.

If $\sigma=G$ is a stable group, then many explicit constructions of the space $X_\tau^G$ exist (see e.g. \cite{rsz}, \cite{ando}, \cite{hosszu} or \cite{sad}); we will now briefly recall one of them: the so-called strengthened Kazarian conjecture proved by Szűcs in \cite{hosszu}. After this we shall construct the classifying space $X_\tau^\i$ of $\tau^\i$-cobordisms.

A stable group $G$ is defined as the direct limit of a sequence of groups $G(n)$ where $G(n)$ has a fixed linear action on $\R^n$ (see definition \ref{stabgr}). We denote by $\gamma_n^G$ the universal rank-$n$ vector bundle with structure group $G$, i.e.
$$\gamma_n^G:=EG(n)\utimes{G(n)}\R^n$$
where $EG(n)\xra{G(n)}BG(n)$ is the universal principal $G(n)$-bundle.

\begin{defi}\label{kazsp}
As before, let $\tau$ be a set of stable singularities and $G$ a stable group. Consider the jet bundle $J_0^\infty(\varepsilon^n,\gamma_{n+k}^G)$ over $BG(n+k)$ which has as fibre the infinite jet space $J_0^\infty(\R^n,\R^{n+k})$ (i.e. the space of all polynomial maps $\R^n\to\R^{n+k}$ with $0$ constant term). Denote by $V_\tau(n)\subset J_0^\infty(\R^n,\R^{n+k})$ the subspace of those maps whose singularity at $0$ is in $\tau$ and let $K_\tau^G(n)$ be the union of the $V_\tau(n)$'s in each fibre of $J_0^\infty(\varepsilon^n,\gamma_{n+k}^G)$, that is,
$$K_\tau^G(n):=EG(n+k)\utimes{G(n+k)}V_\tau(n)$$
where the action of $G(n+k)$ on $V_\tau(n)$ is one-sided, we do not act on the source space $\R^n$ of the polynomial maps. Now the natural map $G(n+k)\to G(n+k+1)$ gives us a map $K_\tau^G(n)\to K_\tau^G(n+1)$. The {\it Kazarian space} for $\tau^G$-maps is the direct limit
$$K_\tau^G:=\lim_{n\to\infty}K_\tau^G(n).$$
\end{defi}

Note that we have a fibration $K_\tau^G(n)\xra{V_\tau(n)}BG(n+k)$.

\begin{defi}
\label{univnb}
The pullback of the vector bundle $\gamma_{n+k}^G\to BG(n+k)$ by the above fibration will be denoted by $\nu_\tau^G(n)\oplus\varepsilon^n\to K_\tau^G(n)$ where $\nu_\tau^G(n)$ is a rank-$k$ virtual vector bundle. Since the natural map $BG(n+k)\to BG(n+k+1)$ induces $\gamma_{n+k}^G\oplus\varepsilon^1$ from $\gamma_{n+k+1}^G$, we get that the map $K_\tau^G(n)\to K_\tau^G(n+1)$ induces $\nu_\tau^G(n)$ from $\nu_\tau^G(n+1)$. The {\it universal virtual normal bundle} of $\tau^G$-maps is
$$\nu_\tau^G:=\lim_{n\to\infty}\nu_\tau^G(n)\to K_\tau^G.$$
\end{defi}

The key tool to proving our main theorems will be the construction of the classifying space of cobordisms of singular maps using the Thom space of the universal virtual normal bundle (for $\tau^G$-cobordisms this is the strengthened Kazarian conjecture \cite[corollary 74]{hosszu} and for $\tau^\i$- and $\tau^{G\oplus\xi}$-cobordisms it will follow from our constructions). Now the Thom space of a virtual vector bundle does not exist, but we can still get a well-defined space (up to homotopy) if we apply the infinite loop space of infinite suspension functor $\Gamma:=\Omega^\infty S^\infty$ of Barratt and Eccles \cite{be} to it. This is a consequence of \cite[remark 60 and definition 72]{hosszu} but it is worthwhile to recall its direct construction here.

\begin{defi}
Let $\nu$ be a virtual vector bundle of rank $k$ over a space $K$ such that $K$ is the direct limit of a sequence of inclusions $K(n)\subset K(n+1)$ where for each $n$ the restriction $\nu(n):=\nu|_{K(n)}$ can be represented by $\xi(n)\ominus\varepsilon^n$ for a (non-virtual) vector bundle $\xi(n)$, moreover, the inclusion $K(n)\into K(n+1)$ is a $c(n)$-homotopy equivalence for a weakly increasing sequence $c\colon\N\to\N$ tending to $\infty$ (this holds for $\nu_\tau^G$ but also more generally for any virtual bundle over a CW complex). Then we can define $S^nT\nu(n)$ as the Thom space $T\xi(n)$, hence $\Gamma T\nu(n)=\Omega^n\Gamma S^nT\nu(n)$ also exists. We define the space 
$$\Gamma T\nu:=\lim_{n\to\infty}\Gamma T\nu(n)=\lim_{n\to\infty}\lim_{m\to\infty}\Omega^{n+m}S^{n+m}T\nu(n).$$
\end{defi}

Note that here the inclusion $S^{n+m}T\nu(n)\into S^{n+m}T\nu(n+1)$ induces isomorphism in cohomologies up to the index $c(n)+k+n+m$ for each $m>0$ (by the Thom isomorphism), hence the $(c(n)+k)$-homotopy type of $\Omega^{n+m}S^{n+m}\nu(n)$ coincides with that of $\Omega^{n+m}S^{n+m}T\nu(n+1)$ (by the Whitehead theorem). This means that the $(c(n)+k)$-homotopy type of $\Omega^{n+m}S^{n+m}T\nu$ can be defined as that of $\Omega^{n+m}S^{n+m}T\nu(n)$. Moreover, we also get that the direct limits $\Gamma T\nu(n)$ and $\Gamma T\nu(n+1)$ are $(c(n)+k)$-homotopy equivalent and since $c(n)$ tends to $\infty$ this implies that for any $r$ the $r$-homotopy type of $\Gamma T\nu$ is that of $\Omega^{n+m}S^{n+m}T\nu(n)$ for some numbers $n$ and $m$.

\begin{rmk}
\label{spectrum}
The space $\Gamma T\nu$ is an infinite loop space, hence it naturally defines a spectrum $E_*:=\Omega^{\infty-*}S^\infty T\nu=\Gamma S^*T\nu$. What is more, the ``virtual space'' $T\nu$ defines an equivalence class of spectra and this $E_*$ is a representative of it; see \cite[remark 62]{hosszu}. 
\end{rmk}

The strengthened Kazarian conjecture \cite[corollary 74]{hosszu} states:

\begin{thm}\label{kazc}
If $\tau$ is a set of stable singularities and $G$ is a stable group, then the classifying space $X_\tau^G$ of cobordisms of $\tau^G$-maps is homotopy equivalent to $\Gamma T\nu_\tau^G$.
\end{thm}

\begin{rmk}
The cohomology theory mentioned in remark \ref{cohrmk} is defined for $\sigma=G$ by the spectrum $E_*^{\tau^G}:=\Gamma S^*T\nu_\tau^G$.
\end{rmk}

\begin{rmk}\label{classmap}
Later it will be important for us to have a deeper understanding of the connection between the cobordism group $\Cob_\tau^G(P)$ (for a manifold $P^{n+k}$) and the classifying homotopy classes $[\cpt P,\Gamma T\nu_\tau^G]$, so we now give a short description of it based on \cite{hosszu}.

If $f\colon M^n\to P^{n+k}$ is a $\tau^G$-map and $i\colon M^n\into\R^r$ is an embedding, then $f\times i\colon M\into P\times\R^r$ is a so-called $\tau^G$-embedding. A $\tau^G$-embedding of $M^n$ into a manifold $Q^{n+k+r}$ is defined as a triple $(e,\VV,\FF)$ where $e\colon M\into Q$ is an embedding with normal $G$-structure, $\VV$ is a sequence $(v_1,\ldots,v_r)$ where the $v_i$'s are pointwise independent vector fields along $e(M)$ (i.e. sections of the bundle $TQ|_{e(M)}$), $\FF$ is a foliation of dimension $r$ on a neighbourhood of $e(M)$ which is tangent to $\VV$ along $e(M)$ and any point $p\in M$ has a neighbourhood on which the composition of $e$ with the projection along the leaves of $\FF$ to a small $(n+k)$-dimensional transverse slice has at $p$ a singularity which belongs to $\tau$.

Cobordisms of $\tau^G$-embeddings of $n$-manifolds to $Q$ can be defined in the usual way and their cobordism group is denoted by $\Emb_\tau^G(n,Q)$. Now if the number $r$ is sufficiently large, then assigning to the $\tau^G$-map $f\colon M^n\to P^{n+k}$ the $\tau^G$-embedding $f\times i\colon M\into P\times\R^r$ (with vector fields arising from a basis of $\R^r$ and foliation composed of the leaves $\{p\}\times\R^r$) yields an isomorphism
$$\Cob_\tau^G(P)\cong\Emb_\tau^G(n,P\times\R^r).$$
This is stated in \cite{hosszu} as theorem 2 and its proof relies on lemma 43 and theorem 1. This is important to note since the proofs of these all only depend on (ambient) isotopies of the embeddings into $P\times\R^r$ which leave the stable normal bundles of the maps involved unchanged.

Now forgetting the singularity structure (given by $\VV$ and $\FF$) of $f\times i$ yields a well-defined cobordism class of an embedding of $M$ into $P\times\R^r$ with normal $G$-structure. Since these are classified by the Thom space $T\gamma_{k+r}^G$ we have the homotopy class of a map $S^r\cpt P\to T\gamma_{k+r}^G$ corresponding to it such that $f\times i(M)$ is the preimage of $BG(k+r)$ and its normal bundle is induced from $\gamma_{r+k}^G$ (note that $S^r\cpt P$ is the one-point compactification of $P\times\R^r$). Then using the singularity structure we can lift this classifying map through the map $T(\nu_\tau^G(r)\oplus\varepsilon^r)\to T\gamma_{k+r}^G$ induced by the fibration $K_\tau^G(r)\to BG(k+r)$ (see definition \ref{univnb}) to get a map
$$\kappa_f\colon S^r\cpt P\to T(\nu_\tau^G(r)\oplus\varepsilon^r)\cong S^rT\nu_\tau^G(r)$$
where $f\times i(M)$ is now the preimage of $K_{\tau}^G(r)$ and its normal bundle is induced from $\nu_\tau^G(r)\oplus\varepsilon^r$.

The adjoint correspondence identifies $[S^r\cpt P,S^rT\nu_{\tau}^G(r)]$ with $[\cpt P,\Omega^rS^rT\nu_{\tau}^G(r)]$ and if $r$ is sufficiently large, the latter is identified with $[\cpt P,\Omega^\infty S^\infty T\nu_{\tau}^G(r)]$ since by increasing $r$ we only attach large dimensional cells to the target space. Now the homotopy class $[\kappa_f]$ can be thought of as an element in $[\cpt P,\Gamma T\nu_\tau^G]$ and assigning $[\kappa_f]$ to $[f]$ gives an isomorphism
$$\Emb_\tau^G(n,P\times\R^r)\cong[\cpt P,\Gamma T\nu_\tau^G].$$

The composition of the two isomorphisms described above show a bijective correspondence between the cobordism classes of the $f$'s and the homotopy classes of the $\kappa_f$'s which forms the isomorphism $\Cob_\tau^G(P)\cong[\cpt P,\Gamma T\nu_\tau^G]$ of theorem \ref{kazc}.
\end{rmk}

Next we shall define the Kazarian space of Wall $\tau$-maps and describe the classifying space $X_\tau^\i$ similarly to the theorem above. Consider the first Stiefel--Whitney class $w_1\colon B\O\to\RP^\infty$ (which is a fibration with fibre $B\SO$) and let $i\colon S^1\into\RP^\infty$ be the inclusion of the $1$-cell.

\begin{defi}\label{bidef}
Let $B\i$ be the homotopy pullback of the diagram $B\O\xra{w_1}\RP^\infty\xleftarrow iS^1$, i.e. the space we obtain by pulling back the fibration $w_1\colon B\O\xra{B\SO}\RP^\infty$ by $i$, and let $w\colon B\i\to S^1$ be the map indicated on the pullback diagram
$$\xymatrix{
B\i\ar[r]^w\ar@{^(->}[d] & S^1\ar@{^(->}[d]^i \\
B\O\ar[r]^{w_1} & \RP^\infty 
}$$
In other words we define $B\i$ as the subspace $S^1\utimes{\Z_2}B\SO$ in $B\O=S^\infty\utimes{\Z_2}B\SO$.
\end{defi}

This $B\i$ is the direct limit of the spaces $B\i(n)$ that we get in a similar way from $B\O(n)$ instead of $B\O$. We can pull back the tautological bundle $\gamma_n^\O\to B\O(n)$ to a vector bundle $\gamma_n^\i\to B\i(n)$ and this makes $B\i(n)$ the classifying space of rank-$n$ vector bundles equipped with Wall structures, that is, those vector bundles $\xi$ for which $w_1(\xi)$ is the mod $2$ reduction of a specific integer cohomology class (this is true since $S^1\cong K(\Z,1)$, $\RP^\infty\cong K(\Z_2,1)$ and $i$ corresponds to the mod $2$ reduction). Now we can again consider the jet bundle $J_0^\infty(\varepsilon^n,\gamma_{n+k}^\i)$ over $B\i(n+k)$ and take in each fibre the space $V_\tau(n)$ to obtain a fibration $K_\tau^\i(n)\xra{V_\tau(n)}B\i(n+k)$ as in definition \ref{kazsp}.

\begin{defi}
If $\tau$ is a set of stable singularities, the {\it Kazarian space} for $\tau^\i$-maps is $K_\tau^\i:=\underset{n\to\infty}\lim K_\tau^\i(n)$.
\end{defi}

We have now a diagram
$$\xymatrix@R=1pc@C=1pc{
K_\tau^\i(n)\ar[rr]\ar[dddr]_{V_\tau(n)}\ar@{^(->}[dr] && K_\tau^\O(n)\ar[dddr]_{V_\tau(n)}\ar@{^(->}[dr] & \\
& J^\infty_0(\varepsilon^n,\gamma^\i_{n+k})\ar[rr]\ar[dd]^{J^\infty_0(\R^n,\R^{n+k})} && J^\infty_0(\varepsilon^n,\gamma^\O_{n+k})\ar[dd]^{J^\infty_0(\R^n,\R^{n+k})} \\
&&&\\
& B\i(n+k)\ar[rr] && B\O(n+k)
}$$
where all squares are pullback squares. Using this we can define the rank-$k$ virtual vector bundle $\nu_\tau^\i(n)\to K_\tau^\i(n)$ as the pullback of $\nu_\tau^\O(n)\to K_\tau^\O(n)$.

\begin{defi}
The {\it universal virtual normal bundle} of $\tau^\i$-maps is $\nu_\tau^\i:=\underset{n\to\infty}\lim\nu_\tau^\i(n)$ over $K_\tau^\i$.
\end{defi}

\begin{thm}\label{tauithm}
If $\tau$ is a set of stable singularities, then we have $X_\tau^\i\cong\Gamma T\nu_\tau^\i$.
\end{thm}

\begin{prf}
  We need an isomorphism $\Cob_\tau^\i(P)\cong[\cpt P,\Gamma T\nu_\tau^\i]$ for any target $P$. Recall from remark \ref{classmap} that the analogous isomorphism for $\tau^\O$-cobordisms is obtained by assigning $\tau$-embeddings to the $\tau$-maps, then lifting their classifying maps through the map $S^rT\nu_\tau^\O(r)\to T\gamma_{r+k}^\O$ (where $r$ is a large number). As we noted in remark \ref{classmap} the isomorphism of the cobordism groups of $\tau$-maps and $\tau$-embeddings (see \cite[theorem 2]{hosszu}) does not depend on the choice of normal structures, hence it also yields the isomorphism
  $$\Cob_\tau^\i(P)\cong\Emb_\tau^\i(n,P\times\R^r)$$
for $r$ sufficiently large, where $\Emb_\tau^\i(n,P\times\R^r)$ is defined analogously.

Now if we have an embedding with normal Wall structure, then its classifying map can be lifted to $T\gamma_{r+k}^\i$ so if it is also a $\tau$-embedding, then using the pullback square
$$\xymatrix{
  S^rT\nu_\tau^\i(r)\ar[r]\ar[d] & S^rT\nu_\tau^\O(r)\ar[d] \\
  T\gamma_{r+k}^\i\ar[r] & T\gamma_{r+k}^\O
}$$
this inducing map lifts to $S^rT\nu_\tau^\i(r)$. So $\tau^\i$-embeddings can be induced from $S^rT\nu_\tau^\i(r)$ and since all of the above extnds to maps between manifolds with boundary (in particular to cobordisms), we get a homomorphism from $\Cob_\tau^\i(P)$ to $[\cpt P,\Gamma T\nu_\tau^\i]$ as in remark \ref{classmap}.

To see that this is an isomorphism, note that if we have a map $S^r\cpt P\to S^rT\nu_\tau^\i(r)$, then its composition with the map to $S^rT\nu_\tau^\O(r)$ pulls back a $\tau$-embedding into $P\times\R^r$ (by theorem \ref{kazc}) which also has a normal Wall structure. Similarly homotopies between maps to $S^rT\nu_\tau^\i(r)$ pull back cobordisms of $\tau$-embeddings with normal Wall structures. Thus we obtained the inverse of the above homomorphism $\Cob_\tau^\i(P)\to[\cpt P,\Gamma T\nu_\tau^\i]$, meaning that it is iso and this concludes our proof.
\end{prf}

\section{A long exact sequence}\label{lessect}

We shall now connect the classifying spaces of $\tau^\i$-cobordisms and $\tau^\SO$-cobordisms. 
Consider the embedding $B\SO\into B\O$ as a fibre over a point in $S^1\subset\RP^\infty$ which then factors through $v\colon B\SO\into B\i$ (again as the embedding of a fibre):
$$\xymatrix@C=.75pc{
&& B\i\ar[d]\ar@{=}[r] & S^1\utimes{\Z_2}B\SO\ar[rr]^(.6)w && S^1\ar@{^(->}[d]^i \\
B\SO\ar@{^(->}[rr]\ar[urr]^v
 && B\O\ar@{=}[r] & S^\infty\utimes{\Z_2}B\SO\ar[rr]^(.6){w_1} && \RP^\infty
}$$

If we now take the approximations $B\SO(n)$, $B\i(n)$ and $B\O(n)$ of the spaces above and the approximation $v(n)\colon B\SO(n)\to B\i(n)$ of $v$, then both squares in the diagram below will be pullback squares:
$$\xymatrix{
\gamma^\SO_n\ar[r]\ar[d] & \gamma^\i_n\ar[r]\ar[d] & \gamma^\O_n\ar[d] \\
B\SO(n)\ar[r]^{v(n)} & B\i(n)\ar[r] & B\O(n)
}$$
But then the same is true for the jet bundles over the spaces $B\SO(n+k)$, $B\i(n+k)$ and $B\O(n+k)$, hence also for the Kazarian spaces and so the suspensions of the approximations of the universal virtual normal bundles also fit into a diagram
$$\xymatrix@R=1pc@C=.7pc{
\nu_\tau^\SO(n)\oplus\varepsilon^n\ar[rr]\ar[dd]\ar[dr] && \nu_\tau^\i(n)\oplus\varepsilon^n\ar[rr]\ar[dd]\ar[dr] && \nu_\tau^\O(n)\oplus\varepsilon^n\ar[dd]\ar[dr] & \\
& \gamma^\SO_{n+k}\ar[rr]\ar[dd] && \gamma^\i_{n+k}\ar[rr]\ar[dd] && \gamma^\O_{n+k}\ar[dd] \\
K_\tau^\SO(n)\ar[rr]\ar[dr] && K_\tau^\i(n)\ar[rr]\ar[dr] && K_\tau^\O(n)\ar[dr] & \\
& B\SO(n+k)\ar[rr]^{v(n+k)} && B\i(n+k)\ar[rr] && B\O(n+k)
}$$
with all squares being pullback squares. After these preliminary observations we are ready to state the following (which is the main statement of theorem \ref{t1}):

\begin{thm}\label{les}
For any set $\tau$ of stable singularities and any manifold $Q^q$ there is a long exact sequence
\begin{alignat*}2
\ldots&\xra{\psi_{m+1}}\Cob_\tau^\SO(Q\times\R^m)\xra{\varphi_m}\Cob^\SO_\tau(Q\times\R^m)\xra{\chi_m}\Cob_\tau^\i(Q\times\R^m)\xra{\psi_m} \\
&\xra{\psi_m}\Cob_\tau^\SO(Q\times\R^{m-1})\to\ldots
\end{alignat*}
\end{thm}

\begin{prf}
Let $A\subset S^1$ be any contractible subspace and let $B:=w^{-1}(A)$ be its preimage in $B\i$. If $v\colon B\SO\into B\i$ was the embedding of the fibre over a point in $A$, then $v(B\SO)\subset B$ and the map $v$ is a homotopy equivalence $B\SO\cong B$. Moreover, if $K$ is the restriction of the Kazarian space $K_\tau^\i$ over $B\subset B\i$, then we have a pullback square
$$\xymatrix{
K_\tau^\SO\ar[d]\ar[r] & K_\tau^\i|_B=K\ar[d] \\
B\SO\ar[r]^\cong & B
}$$
with the vertical arrows being fibrations, hence we have $K\cong K_\tau^\SO$. Similarly if $\nu$ is the restriction of $\nu_\tau^\i$ over $K\subset K_\tau^\i$, then we have a pullback square 
$$\xymatrix{
\nu_\tau^\SO\ar[d]\ar[r] & \nu_\tau^\i|_K=\nu\ar[d] \\
K_\tau^\SO\ar[r]^\cong & K
}$$
which implies that $\nu$ is stably isomorphic to $\nu_\tau^\SO$.

Now let $A_1$ be any point $p\in S^1$ and put $A_2:=S^1\setminus\{p\}$ (which are contractible subspaces). Let the preimages of $A_1$ and $A_2$ under $w$ respectively be $B_1$ and $B_2$ (thus we have $B_1\cong B\SO\cong B_2$), their preimages in $K_\tau^\i$ be $K_1$ and $K_2$ (thus we have $K_1\cong K_\tau^\SO\cong K_2$) and put $\nu_1:=\nu_\tau^\i|_{K_1}$ and $\nu_2:=\nu_\tau^\i|_{K_2}$ (which are now stably isomorphic to $\nu_\tau^\SO$). Then for all $n$ we have a cofibraion
$$S^nT\nu_2(n)\into S^nT\nu_\tau^\i(n)\to S^nT\nu_\tau^\i(n)/S^nT\nu_2(n)$$
(where $\nu_1(n)$ and $\nu_2(n)$ are the appropriate restrictions of $\nu_\tau^\i(n)$). Since the normal bundle of $K_1\subset K_\tau^\i$ is induced from the normal bundle of $p\in S^1$ it is trivial and so we have $S^nT\nu_\tau^\i(n)/S^nT\nu_2(n)=S^{n+1}T\nu_1(n)$, hence the cofibration above has the form
$$S^nT\nu_\tau^\SO(n)\into S^nT\nu_\tau^\i(n)\to S^{n+1}T\nu_\tau^\SO(n).$$

Now applying the functor $\Omega^{n+m}\Gamma$ to the Puppe sequence of this cofibration we get a sequence of maps
$$\ldots\to\Omega^m\Gamma T\nu_\tau^\SO(n)\to\Omega^m\Gamma T\nu_\tau^\i(n)\to\Omega^m\Gamma ST\nu_\tau^\SO(n)\to\Omega^{m-1}\Gamma T\nu_\tau^\SO(n)\to\ldots$$
And note that here $\Omega^m\Gamma ST\nu_\tau^\SO(n)$ equals $\Omega^{m-1}\Gamma T\nu_\tau^\SO(n)$. This sequence is infinite to the right by construction, but it is also infinite to the left since the number $n$ could be arbitrary and we get the same maps by applying $\Omega^{n+m}\Gamma$ to the $n$'th suspensions as by applying $\Omega^{n+m+1}\Gamma$ to the $(n+1)$'st suspensions. We also note that the maps in this sequence commute with the natural maps $\Omega^m\Gamma T\nu_\tau^\sigma(n)\to\Omega^m\Gamma T\nu_\tau^\sigma(n+1)$ (for $\sigma=\SO,\i$). Then converging with $n$ to infinity yields a sequence
$$\ldots\to\Omega^m\Gamma T\nu_\tau^\SO\to\Omega^m\Gamma T\nu_\tau^\i\to\Omega^{m-1}\Gamma T\nu_\tau^\SO\to\Omega^{m-1}\Gamma T\nu_\tau^\SO\to\ldots$$
of the direct limits.

If we then fix a manifold $Q$ and apply the functor $[\cpt Q,\cdot]$ to this sequence, then we obtain the long exact sequence of cobordism groups as claimed.
\end{prf}

\begin{rmk}
\label{cofib}
In the proof above we used the cofibrations
$$S^nT\nu_\tau^\SO(n)\into S^nT\nu_\tau^\i(n)\to S^{n+1}T\nu_\tau^\SO(n)$$
for all $n$. Later we will refer to this (and to similar occurences) as a cofibration
$$T\nu_\tau^\SO\into T\nu_\tau^\i\to ST\nu_\tau^\SO.$$
Although these are not actual spaces, applying the functor $\Gamma$ to them yields a well-defined fibration (cf. the proof of \cite[theorem 8]{hosszu}) which is a fibration between the corresponding spectra described in remark \ref{spectrum}.
\end{rmk}

\begin{rmk}
If we had in theorem \ref{les} $Q=\R^q$, then the same long exact sequence could be obtained from the cofibration
$$T\nu_\tau^\SO\into T\nu_\tau^\i\to ST\nu_\tau^\SO$$
by just applying the functor $\Gamma$ which yields a fibration
$$\Gamma T\nu_\tau^\i\xra{\Gamma T\nu_\tau^\SO}\Gamma ST\nu_{\tau}^\SO$$
and then taking the homotopy long exact sequence of this.
\end{rmk}

Now to prove theorem \ref{t1} we only have to describe the homomorphisms in the above exact sequence.

\section{Description of the homomorphisms $\varphi_m$, $\chi_m$ and $\psi_m$}\label{phi}

First fix a manifold $K$, a vector bundle $\zeta\to K$ of rank $r$, a closed submanifold $A\subset K$ of codimension $m$ and with normal bundle $\xi$ and put $B:=K\setminus A$. If we denote by $\Emb^\zeta(n,P^{n+r})$ the cobordism semigroup of $n$-manifolds embedded into a fixed manifold $P^{n+r}$ with normal bundle induced from $\zeta$ (i.e. $\zeta$-embeddings), then we have
$$\Emb^\zeta(n,P)\cong[\cpt P,T\zeta].$$
Hence the Puppe sequence of the cofibration $T\zeta|_B\into T\zeta\to T\zeta/T\zeta|_B=T(\zeta|_A\oplus\xi)$ gives us a sequence of homomorphisms
$$\Emb^{\zeta|_B}(n,P)\to\Emb^\zeta(n,P)\to\Emb^{\zeta|_A\oplus\xi}(n-m,P)\xra\partial\Emb^{\zeta|_B\oplus\varepsilon^1}(n-1,P)\to\ldots$$
of cobordism semigroups (which is exact in the stable dimensions).

\begin{lemma}\label{iotal}
Let $f\colon M^{n-m}\into P^{n+r}$ be a $(\zeta|_A\oplus\xi)$-embedding which yields a decomposition $\nu_f=(\zeta|_A)_M\oplus\xi_M$ (with $(\zeta|_A)_M$ and $\xi_M$ induced from $\zeta|_A$ and $\xi$ respectively) and the exponential map can be used to extend $f$ to an embedding $\tilde f$ of the total space of $\xi_M$ into a small neighbourhood of $f(M)$. Then the boundary homomorphism $\partial$ assigns to the cobordism class $[f]$ the class of $\tilde f|_{S\xi_M}$ which has codimension $r+1$, is endowed with the natural normal vector field of the sphere bundle and the orthogonal complement of this trivial subbundle in its normal bundle is induced from $\zeta|_B$.
\end{lemma}

\begin{prf}
We obtain $T(\zeta|_A\oplus\xi)$ in the Puppe sequence by restricting $\zeta$ to a tubular neighbourhood $U\subset K$ of $A$ and then factoring $T\zeta$ by the complement of this restriction, that is, attaching the cone $CT(\zeta|_{K\setminus U})$. The next item in the sequence is $ST(\zeta|_B)$ which we get by gluing the cone $CT\zeta$ to the previously obtained space. In the following we shall always mean the above representations of the spaces $T(\zeta|_A\oplus\xi)$ and $ST(\zeta|_B)$.

\begin{figure}[H]
\begin{center}
\centering\includegraphics[scale=0.1]{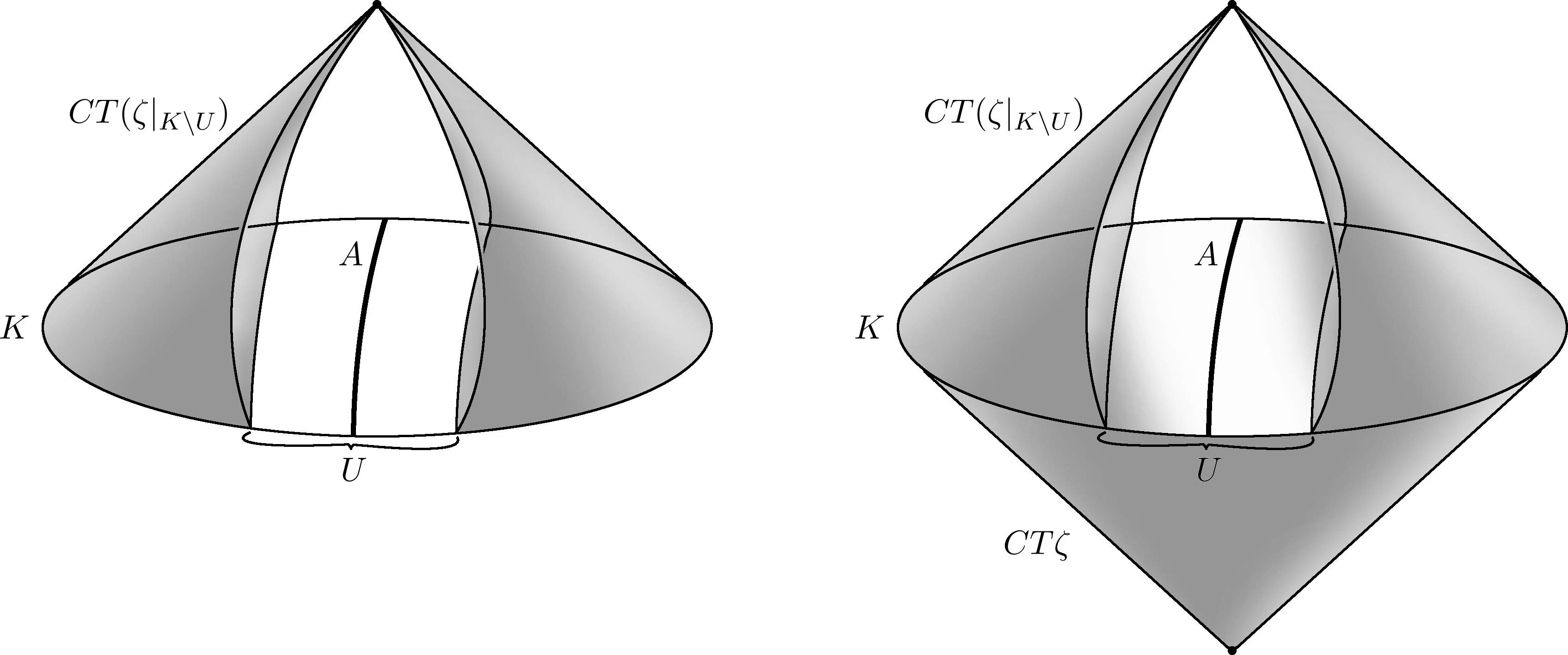}
\begin{changemargin}{2cm}{2cm} 
\caption{\hangindent=1.4cm\small We show the local form of $T(\zeta|_A\oplus\xi)$ (on the left) and $ST(\zeta|_B)$ (on the right) in the case $r=0$, $m=1$.}\label{kep1}
\end{changemargin} 
\vspace{-1.3cm}
\end{center}
\end{figure}

If the homotopy class of $\kappa_f\colon\cpt P\to T(\zeta|_A\oplus\xi)$ induces the cobordism class of $f\colon M\into P$ (such that $f(M)=\kappa_f^{-1}(A)$), then its image under
$$\partial\colon[\cpt P,T(\zeta|_A\oplus\xi)]\to[\cpt P,ST(\zeta|_B)]$$
is represented by $\lambda_f\colon\cpt P\to ST(\zeta|_B)$ which we get by taking the portion of the image of $\kappa_f$ in the total space of $\zeta|_U$ and pushing it toward the special point of the cone $CT(\zeta|_U)\subset CT\zeta$, in this way detaching it from $A$, and otherwise leaving the rest of $\kappa_f$ unchanged.

\begin{figure}[H]
\begin{center}
\centering\includegraphics[scale=0.1]{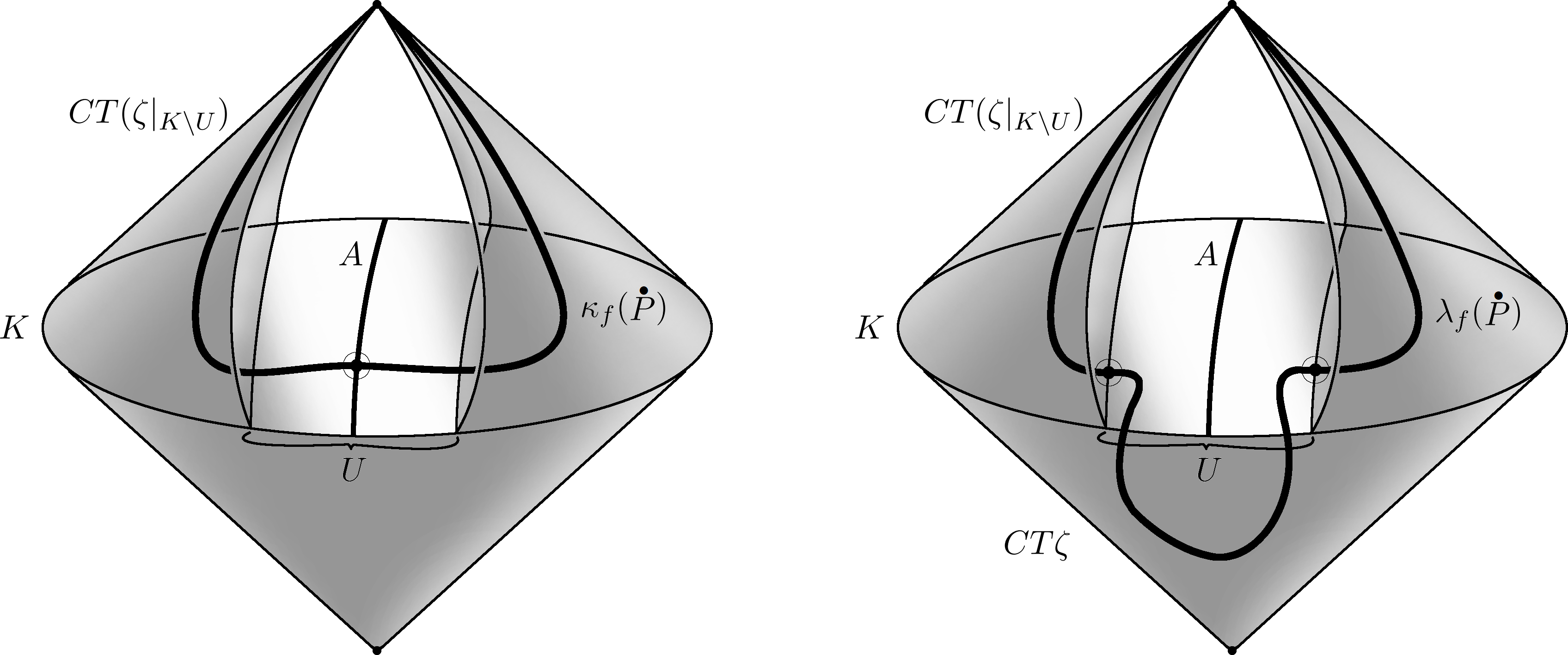}
\begin{changemargin}{2cm}{2cm} 
\caption{\hangindent=1.4cm\small We show the local form of the image of $\kappa_f$ (on the left) and $\lambda_f$ (on the right) in $ST(\zeta|_B)$ again in the case $r=0$, $m=1$.}\label{kep2}
\end{changemargin} 
\vspace{-1.3cm}
\end{center}
\end{figure}

The inducing map $\kappa_f$ is such that the exponential image of the normal bundle $\nu_f=(\zeta|_A)_M\oplus\xi_M$ (restricted to a small neighbourhood) is the preimage of a tubular neighbourhood of $A\subset T(\zeta|_A\oplus\xi)$ and it is mapped fibrewise bijectively to its image in this neighbourhood. Thus we may assume that the image of $\tilde f$ coincides with $\kappa_f^{-1}(U)$ such that $\tilde f(S\xi_M)=\kappa_f^{-1}(\partial U)$ and $\kappa_f$ is again a fibrewise bijection from $S\xi_M$ to $\partial U|_{\kappa_f(\kappa_f^{-1}(A))}$.

Then we can assume that the image of $\lambda_f$ intersects the base space $K$ only in $\partial U$ (this intersection is in $B=K\setminus A$) and $\lambda_f$ restricted to the preimage of this intersection coincides with $\kappa_f$. Hence for any point $p\in(\kappa_f\circ f)^{-1}(A)$ we have that $\lambda_f$ maps the sphere $\tilde f(S_p\xi_M)$ to $B$ and the union of these (that is, $\tilde f(S\xi_M)$) is the whole $\lambda_f^{-1}(B)$. Since $S\xi_M$ is a sphere bundle, its embedding has a natural normal vector field. This trivialises the subbundle of its normal bundle induced from the direction of the suspension in $ST(\zeta|_B)$. The orthogonal complement of this in the normal bundle of $\tilde f(S\xi_M)$ is then induced from the $\zeta|_B$ part of the normal bundle of $B\subset ST(\zeta|_B)$. This is what we wanted to prove.
\end{prf}\medskip

If the manifold above has the form $K=S^1\utimes{\Z_2}V$ and the submanifold $A$ is a fibre $V$, then we have $\xi=\varepsilon^1$, $B\cong V$ and $\zeta|_B\cong\zeta|_A=\zeta|_V$, hence the sequence above has the form
$$\Emb^{\zeta|_V}(n,P)\to\Emb^\zeta(n,P)\to\Emb^{\zeta|_V\oplus\varepsilon^1}(n-1,P)\xra\partial\Emb^{\zeta|_V\oplus\varepsilon^1}(n-1,P)\to\ldots$$

\begin{crly}\label{iotac}
The boundary homomorphism $\partial$ in this case is $\id+\iota$ where $\iota$ is the involution which acts on a cobordism class represented by a $1$-framed $\zeta|_V$-embedding $f\colon M^{n-1}\into P^{n+r}$ such that it inverts its framing vector field $v$ and in the orthogonal complement of $v$ in $\nu_f$ (induced from $\zeta|_V$) it applies the $\Z_2$-action of $V$.
\end{crly}

\begin{prf}
We want to understand the embedding $\tilde f|_{S\xi_M}$ in the lemma above in this special case. Now $\xi$ is a trivial line bundle, so $S\zeta_M=M_-\sqcup M_+$ where the $M_\pm$ are both diffeomorphic to $M$, the maps $f_\pm:=\tilde f|_{M_\pm}$ can be identified with $f$ and the normal vector fields $v_\pm$ giving the $1$-framing of these two maps are opposite since these are the natural (say outward pointing) vectors of a sphere bundle. Thus we may assume that the $1$-framing $v_+$ of $f_+$ is the same as the $1$-framing $v$ of $f$ and $v_-$ is opposite to it.

What is left is to understand how the orthogonal complement of $v_\pm$ in the normal bundle is induced from $\zeta|_V$ on the two manifolds $M_\pm$. If our fixed fibre was $A=\{p\}\times V$ for a point $p\in S^1$, then (with the notation of the previous proof) the tubular neighbourhood $U\subset K$ is $(p-\varepsilon,p+\varepsilon)\times V$ for a small number $\varepsilon$. Now $f(M)=\kappa_f^{-1}(\{p\}\times V)$ and the part of $\nu_f$ induced from $\zeta|_V$ is pulled back by the restriction of $\kappa_f$ from the bundle over $\{p\}\times V$. Identifying the fibres in $[p-\varepsilon,p+\varepsilon]\times V$ in the trivial way then yields that the $\zeta|_V$ parts of the normal bundles of $f_\pm(M_\pm)=\kappa_f^{-1}(\{p\pm\varepsilon\}\times V)$ are induced by the same map from the bundle over $\{p\pm\varepsilon\}\times V$. But then identifying the fibres over $S^1\setminus\{p\}$ (i.e. in $B$) instead of $[p-\varepsilon,p+\varepsilon]$ changes the identification of $\{p+\varepsilon\}\times V$ with $\{p-\varepsilon\}\times V$ by applying the $\Z_2$-action in one of them, say in $\{p-\varepsilon\}\times V$.

Thus we got that both the $1$-framing of $f_+$ and the inducing map of the rest of its normal bundle coincide with those of $f$, on the other hand, the $1$-framing of $f_-$ is opposite to that of $f$ and the inducing map of the rest of its normal bundle is composed with the $\Z_2$-action in $V$. Hence we have $\partial[(f,v)]=[(f_+,v_+)]+[(f_-,v_-)]$ where $[(f_+,v_+)]=[(f,v)]$ and $[(f_-,v_-)]=\iota[(f,v)]$ for the involution $\iota$ described in our statement above.
\end{prf}

\begin{prop}\label{2}
$\varphi_m\colon\Cob_\tau^\SO(Q\times\R^m)\to\Cob_\tau^\SO(Q\times\R^m)$ is of the form $\id+\iota$ where $\iota$ is the involution which, when representing cobordism classes by $\tau^\SO$-embeddings (see remark \ref{classmap}), composes a map with the reflection to a hypersurface $Q\times\R^{m+r-1}\subset Q\times\R^{m+r}$ but does not change its orientation.
\end{prop}

\begin{prf}
Let $q$ be the dimension of $Q$ and put $n:=q+m-k$. Then remark \ref{classmap} and theorem \ref{tauithm} show that for $\sigma=\i,\SO$ the group $\Cob_\tau^\sigma(Q^q\times\R^m)$ is naturally isomorphic to
$$\Emb_\tau^\sigma(n,Q\times\R^{m+r})\cong\Emb^{\nu_\tau^\sigma(r)\oplus\varepsilon^r}(n,Q\times\R^{m+r})$$
where $r$ is a sufficiently large number. Observe also that we have $K_\tau^\i=S^1\utimes{\Z_2}K_\tau^\SO$ and, using the notation of the proof of theorem \ref{les}, we let $K_1$ be the fibre of $K_\tau^\i$ over a point in $S^1$ and $K_2$ is its complement in $K_\tau^\i$ and for $i=1,2$ we let the virtual bundle $\nu_i$ be the restriction of $\nu_\tau^\i$ over $K_i$ (hence $K_i\cong K_\tau^\SO$ and $\nu_i\cong\nu_\tau^\SO$). We now think of these Kazarian spaces as finite dimensional approximations of the actual spaces (hence closed manifolds) over which the $(r-1)$'st suspensions of the universal virtual normal bundles exist as (non-virtual) vector bundles, but for simplicity of notation we are not indicating this.

Then the source and target of $\varphi_m$ are
$$[S^{m+r}\cpt Q,S^rT\nu_1]\quad\text{and}\quad[S^{m+r}\cpt Q,S^rT\nu_2]$$
respectively, thus we are in the setting of corollary \ref{iotac} with $K=K_\tau^\i$, $V=K_\tau^\SO$, $A=K_1$, $B=K_2$ and $\zeta=\nu_\tau^\i\oplus\varepsilon^r$ with $\zeta|_V=\nu_\tau^\SO\oplus\varepsilon^r$, and so the boundary homomorphism $\partial=\varphi_m$ is indeed $\id+\iota$. The involution $\iota$ here inverts one (say, the last) of the framing vectors of any $\tau^\SO$-embedding into $Q\times\R^{m+r}$ (see remark \ref{classmap}) and applies the $\Z_2$-action of $K_\tau^\SO$ in the orthogonal complement of this vector in the normal bundle.

But if $r$ is large enough, then we may assume for any $\tau^\SO$-embedding to map into a hypesurface $Q\times\R^{r+m-1}\subset Q\times\R^{r+m}$ with the last framing vector being the normal vector of this hypersurface, hence inverting it corresponds to the reflection. This would change the orientation of the normal bundle but the $\Z_2$-action of $K_\tau^\SO$ reverses the orientation in the orthogonal complement too, which means that the orientation of the normal bundle remains unchanged. This finishes the proof.
\end{prf}

\begin{prop}
$\chi_m\colon\Cob_\tau^\SO(Q\times\R^m)\to\Cob_\tau^\i(Q\times\R^m)$ is the forgetful homomorphism that assigns to an oriented cobordism class $[f]$ the class of $(f,0)$ as a Wall map.
\end{prop}

\begin{prf}
This follows immediately since the map between the classifying spaces is just the inclusion of $\Gamma T\nu_2$ into $\Gamma T\nu_\tau^\i$ (see the proof of theorem \ref{les}).
\end{prf}

\begin{prop}\label{pd}
$\psi_m\colon\Cob_\tau^\i(Q\times\R^m)\to\Cob_\tau^\SO(Q\times\R^{m-1})$ assigns to a cobordism class $[(f,w)]$ the oriented cobordism class of $f$ restricted to the Poincaré dual of $w$ which can be represented by a $1$-codimensional submanifold uniquely up to cobordism, the restriction $f|_{\PD(w)}$ can be assumed to be a mapping to $Q\times\R^{m-1}$, it has oriented normal bundle and its cobordism class only depends on the class of $(f,w)$, hence this assignment is well-defined.
\end{prop}

\begin{prf}
The cobordism class $[(f,w)]$ is represented by a $\tau$-map $f\colon M\to Q\times\R^m$ together with a fixed cohomology class $w\in H^1(M;\Z)$ such that its mod $2$ reduction is $w_1(\nu_f)$. Then we have a commutative diagram
$$\xymatrix@C=3.5pc{
& S^1\ar@{^(->}[d] \\
M\ar[ur]^{w}\ar[r]^{w_1(\nu_f)} & \RP^\infty
}$$

The Poincaré dual of 
$w$ is represented by the codimension-$1$ submanifold $N:=w^{-1}(p)\subset M$ for a point $p\in S^1$. Note that the normal bundle of $N$ in $M$ is trivial and so $f|_N$ has a normal vector field which is uniquely defined by a fixed normal vector of $p$ in $S^1$; we also have $w_1(\nu_{f|_N})=w_1(\nu_f)|_N=0$, hence $f|_N$ is orientable and can be canonically oriented (since the classifying map of the stabilisation of $\nu_f$ maps to $B\i$, hence that of $\nu_{f|_N}$ maps to the fibre $B_1\cong B\SO$ of $B\i$ over $p\in S^1$); finally we may assume that $N$ intersects each singularity stratum in $M$ transversally, thus $f|_N$ is a $\tau^{\SO\oplus1}$-map.

By (\ref{+m}) we have $\Cob_\tau^{\SO\oplus1}(Q\times\R^m)\cong\Cob_\tau^\SO(Q\times\R^{m-1})$ which means that assigning to the cobordism class of $(f,w)$ the class of $f|_N$ gives a map
$$\Cob_\tau^\i(Q\times\R^m)\to\Cob_\tau^\SO(Q\times\R^{m-1})$$
taking $[f]\mapsto[f|_{\PD(w)}]$. We now only need to prove that this map is well-defined and is the same as $\psi_m$ for which it is sufficient to prove that $\psi_m$ assigns to the classifying map of $f$ the classifying map of $f|_N$.

We again refer to remark \ref{classmap}: $[f]$ bijectively corresponds to a homotopy class $[\kappa_f]\in[S^m\cpt Q,\Gamma T\nu_\tau^\i]$ represented by $\kappa_f\colon S^{r+m}\cpt Q\to S^rT\nu_\tau^\i$ where $r$ is a large number and we again consider such a finite dimensional approximation of the Kazarian space over which $S^rT\nu_\tau^\i$ exists. In the proof of theorem \ref{les} we obtained the map in the Puppe sequence inducing $\psi_m$ by factoring $T(\nu_\tau^\i\oplus\varepsilon^r)$ by its restriction over the complement of a tubular neighbourhood of $K_1\subset K_\tau^\i$. Recall that $K_1\cong K_\tau^\SO$ is the preimage of $B_1\cong B\SO\subset B\i$ (the fibre of $B\i$ over $p\in S^1$). Thus we have a diagram
$$\xymatrix@R=1pc@C=1pc{
T(\nu_\tau^\i\oplus\varepsilon^r)\ar[rr]\ar[dr]^q && T\gamma^\i_{r+k}\ar[dr] &&& \\
& ST(\nu_1\oplus\varepsilon^r)\ar[rr] && ST\gamma^\SO_{r+k} && \\
K_\tau^\i\ar[rr]\ar@{^(->}[uu] && B\i(r+k)\ar[rr]\ar@{^(->}[uu] && S^1 & \\
& K_1\ar@{_(->}[ul]\ar[rr]\ar@{^(->}[uu] && B_1(r+k)\ar@{_(->}[ul]\ar[rr]\ar@{^(->}[uu] && p\ar@{{}{}{}}[ul]|{\text{\rotatebox{130}{$\in$}}}
}$$
where the top two rows are connected by quotient maps which factor by the restrictions of the Thom spaces over complements of tubular neighbourhoods.

The homomorphism $\psi_m$ assigns to the homotopy class of the inducing map $\kappa_f$ the homotopy class of its composition with the quotient map $q$. Now $\kappa_f$ pulls back the $\tau^\i$-embedding $f\times i$ (for an embedding $i\colon M\into\R^r$) from the embedding of $K_\tau^\i$ into $T(\nu_\tau^\i\oplus\varepsilon^r)$ and the composition $q\circ\kappa_f$ pulls back a representative of $\psi_m[f\times i]$ from the embedding of $K_1$ into $ST(\nu_1\oplus\varepsilon^r)$. But note that we have $N=(q\circ\kappa_f)^{-1}(K_1)$ and its $\tau^\SO$-embedding into $Q\times\R^{r+m}$ is $f|_N\times i$. Hence we proved that $\psi_m[f]=[f|_N]$ and this is what we needed.
\end{prf}\medskip

This finishes the proof of theorem \ref{t1}.

\begin{rmk}\label{gen1}
If the codimension $k$ is sufficiently large (compared to $n$), then we have $\Cob_\tau^\SO(\R^{n+k})=\Omega_n$ and $\Cob_\tau^\i(\R^{n+k})=\WW_n$ for any (non-empty) singularity set $\tau$. Hence in the case $Q\times\R^m=\R^{n+k}$ the portion of the long exact sequence in theorem \ref{t1} where $n$ is sufficiently small (compared to $k$) looks like
$$\ldots\to\Omega_n\to\Omega_n\to\WW_n\to\Omega_{n-1}\to\ldots$$
If we fix here $n$, then $k$ can be increased arbitrarily and the homomorphisms in this sequence do not change, hence this sequence is infinite both to the right and to the left. Moreover, by the above propositions these homomorphisms can be identified with those in the classical exact sequence (\ref{ces1}), hence theorem \ref{t1} generalises the sequence (\ref{ces1}).
\end{rmk}

\begin{rmk}
It is easy to see that propositions \ref{2} and \ref{pd} also extend \cite[lemma 2]{li}.
\end{rmk}

\part{The unoriented case}\label{p2}

\section{Cobordism of $\tau^{G\oplus\xi}$-maps}\label{tauxi}

Fix a stable group $G$ which is the direct limit of the groups $G(n)$. Recall that $G(n)$ has a fixed linear action on $\R^n$ and $\gamma_n^G\to BG(n)$ denotes the universal rank-$n$ vector bundle with structure group $G$. Fix also a rank-$m$ vector bundle $\xi$ over the space $BG=\underset{n\to\infty}\lim BG(n)$.

\begin{defi}\label{deftauxi}
Let $\tau$ be a set of stable singularities. We define a {\it $\tau^{G\oplus\xi}$-map} from a compact manifold $M^n$ to a manifold $P^{n+k+m}$ to be the germ along the zero-section of a $\tau^G$-map
$$\tilde{f}\colon\xi_M\to P$$
where $\xi_M$ is a rank-$m$ vector bundle over $M$ and $\tilde f$ has the following properties:
\begin{itemize}
\item[(i)] the differential $d\tilde f$ restricted to any fibre of $\xi_M$ is injective,
\item[(ii)] noting that (i) implies $\nu_{(\tilde f|_M)}=(\nu_{\tilde f})|_M\oplus\xi_M$ we require that 
the map $M\to BG$ inducing the (virtual) $G$-bundle $(\nu_{\tilde f})|_M$ also pulls back $\xi_M$ from $\xi$, thus it pulls back $\nu_{(\tilde f|_M)}\oplus\varepsilon^r$ from $\gamma^G_{k+r}\oplus\xi$ for sufficiently large numbers $r$.
\end{itemize}
The cobordism of two $\tau^{G\oplus\xi}$-maps is defined in the usual way (i.e. as a $\tau^{G\oplus\xi}$-map to $P\times[0,1]$ where both the germ and the map to $BG$ inducing the vector bundle extend those of the boundary) and the cobordism group of $\tau^{G\oplus\xi}$-maps to the manifold $P$ is denoted  by $\Cob_\tau^{G\oplus\xi}(P)$.
\end{defi}


\begin{ex}
In the case $G=\SO$ and $\xi=\varepsilon^m$, the above definition is exactly the definition of $m$-framed $\tau$-maps in \cite{hosszu}.
\end{ex}

As we saw previously in part \ref{p1}, the main tools for understanding how cobordisms of $\tau$-maps are induced from the classifying spaces are $\tau$-embeddings. In order to construct the classifying space of $\tau^{G\oplus\xi}$-maps we now introduce the following:

\begin{defi}
Let $\tau$ be a set of stable singularities. By a {\it $\tau^{G\oplus\xi}$-embedding} of a manifold $M^n$ into a manifold $Q^q$ we mean a quadruple $(e,\VV,\xi_M,\FF)$ where
\begin{itemize}
\item[(i)] $e\colon M\into Q$ is an embedding,
\item[(ii)] $\VV=(v_1,\ldots,v_r)$ where $r=q-n-k-m$ and the $v_i$'s are pointwise independent vector fields along $e(M)$, i.e. sections of the bundle $TQ|_{e(M)}$; we identify $\VV$ with the trivialised subbundle generated by the $v_i$'s,
\item[(iii)] $\xi_M$ is a rank-$m$ subbundle of $TQ|_{e(M)}$ which is pointwise independent of both $\VV$ and $Te(M)$ and the normal bundle $\nu_e=\nu'\oplus\xi_M$ is induced from the bundle $\gamma_{k+r}^G\oplus\xi|_{BG(k+r)}$ over $BG(k+r)$ by a map $M\to BG(k+r)$ which pulls back $\xi_M$ from $\xi|_{BG(k+r)}$ and $\nu'$ from $\gamma_{k+r}^G$,
\item[(iv)] $\FF$ is a foliation of dimension $r+m$ on a neighbourhood of $e(M)$ and it is tangent to $\VV\oplus\xi_M$ along $e(M)$,
\item[(v)] any point $p\in M$ has a neighbourhood on which the composition of $e$ with the projection along the leaves of $\FF$ to a small $(n+k)$-dimensional transverse slice has at $p$ a singularity which belongs to $\tau$.
\end{itemize}
The cobordism of two $\tau^{G\oplus\xi}$-embeddings is defined in the usual way and the cobordism group of $\tau^{G\oplus\xi}$-embeddings of $n$-manifolds to $Q$ is denoted  by $\Emb_\tau^{G\oplus\xi}(n,Q)$.
\end{defi}

\begin{rmk}
Such a $\tau^{G\oplus\xi}$-embedding induces a stratification of $M$ by the submanifolds
$$\eta(e):=\eta(e,\VV,\xi_M,\FF):=\{p\in M\mid p\in\eta(\pi\circ e)\}$$
where the $[\eta]$ are the elements of $\tau$ and $\pi$ denotes the local projection around $e(p)$ along the leaves of $\FF$.
\end{rmk}

\begin{ex}\label{vert}
If $\xi_M$ is a vector bundle over $M^n$, $\tilde f\colon\xi_M\to P^{n+k+m}$ is a $\tau^{G\oplus\xi}$-map and $i\colon M\into\R^r$ is an embedding, then we can define a $\tau^{G\oplus\xi}$-embedding $(e,\VV,\xi_M,\FF)$ of $M$ into $P\times\R^r$: We choose an arbitrarily small representative $f$ of $\tilde f$ and put $e:=f|_M\times i$; the vector fields $v_i$ arise from a basis in $\R^r$; the bundle $\xi_M$ can now be viewed as a subbundle of $T(P\times\R^r)|_{e(M)}$ (since for any $p\in M$ we can identify $\xi_M|_p$ with the translate of $df_p(\xi_M|_p)$ to $T_{e(p)}(P\times\R^r)$); and $\FF$ is composed of the leaves $f(\xi_M|_p)\times\R^r$ intersected by a small neighbourhood of $e(M)$ (for all $p\in M$).
\end{ex}

\begin{defi}
~\vspace{-.5em}

\begin{enumerate}
\item The vector fields $\VV=(v_1,\ldots,v_r)$ and the foliation $\FF$ in the above example are called {\it vertical} in $P\times\R^r$.
\item The subsets $P\times\{x\}\subset P\times\R^r$ (for any $x\in\R^r$) are called {\it horizontal sections}.
\end{enumerate}
\end{defi}

The rest of this section consists of constructing the classifying space of the cobordism group of $\tau^{G\oplus\xi}$-maps in a similar way as sketched in remark \ref{classmap}. 
First we state a lemma and a theorem which connect $\tau^{G\oplus\xi}$-maps with $\tau^{G\oplus\xi}$-embeddings and give a stabilisation property for these. Their proofs are the direct analogues of the proofs of \cite[lemma 43]{hosszu} and \cite[theorem 2]{hosszu}
, we just repeat them here for completeness.

\begin{lemma}
Let $\tau$ be a set of stable singularities and $(e,\VV,\xi_M,\FF)$ be a $\tau^{G\oplus\xi}$-embedding of $M^n$ into $P^{n+k+m}\times\R^r$ where $M$ is a compact manifold and $P$ is any manifold. Then there is a diffeotopy $\Phi_t~(t\in[0,1])$ of $P\times\R^r$ such that $\Phi_0$ is the identity and (the differential of) $\Phi_1$ takes $\VV$ to the vertical vector fields $\VV'$ and $\FF$ to the vertical foliation $\FF'$ around the image of $M$. The relative version of this claim is also true, that is, if the vector fields $\VV$ and the foliation $\FF$ are already vertical on a neighbourhood of a compact subset $C\subset e(M)$, then the diffeotopy $\Phi_t~(t\in[0,1])$ is fixed on a neighbourhood of $C$.
\end{lemma}

\begin{prf}
The manifold $M$ is finitely stratified by the submanifolds
$$S_i:=\bigcup_{[\eta]\in\tau\atop\dim\eta(e)=i}\eta(e),\quad i=0,\ldots,n.$$
By the stratified compression theorem (\cite[theorem 1]{hosszu}, the analogue of the multi-compression theorem of \cite{rs} for stratified manifolds) there is a diffeotopy of $P\times\R^r$ which turns the vector fields $\VV$ into vertical vector fields. Therefore we may assume that $\VV$ is already vertical and so we can also assume that the fibres of $\xi_M$ are horizontal (i.e. they are in the tangent spaces of horizontal sections). Hence we only need to find a diffeotopy that takes the foliation $\FF$ into the vertical foliation $\FF'$ and its differential keeps the vector fields $\VV=\VV'$ vertical.

We will recursively deform $\FF$ into $\FF'$ around the images of the strata $S_i~(i=0,\ldots,n)$. First we list some general observations:
\begin{enumerate}
\item\label{tr1} If $R\subset\R^r$ and $L,L'\subset P\times R$ are such that each of $L$ and $L'$ intersects each horizontal section $P\times\{x\}~(x\in R)$ exactly once, then a bijective correspondence $L\to L'$ arises by associating the points on the same horizontal section to each other.
\item\label{tr2} If $A\subset P\times\R^r$ is such that for each $a\in A$ subsets $R_a,L_a,L'_a$ are given as in \ref{tr1}, then a family of bijective maps $\{L_a\to L'_a\mid a\in A\}$ arises. If we have $L_{a_1}\cap L_{a_2}=\varnothing=L'_{a_1}\cap L'_{a_2}$ for any two different points $a_1,a_2\in A$, then the union of these bijections gives a continuous bijective map
$$\varphi\colon U:=\bigcup_{a\in A}L_a\to\bigcup_{a\in A}L'_a=:U'$$
\item\label{tr3} If the subsets $A,L_a,L'_a$ in \ref{tr2} are submanifolds of $P\times\R^r$ such that $U$ and $U'$ are also submanifolds, then the map $\varphi$ is smooth. In this case for all points $(p,x)\in U$ we can join $(p,x)$ and $\varphi(p,x)$ by a minimal geodesic in the horizontal section $P\times\{x\}$, and using these we can extend $\varphi$ to an isotopy $\varphi_t~(t\in[0,1])$ of $U$ (for which $\varphi_0=\id_U$ and $\varphi_1=\varphi$).
\end{enumerate}

Denote by $\VV^\perp$ the orthogonal complement of the bundle $\VV|_{e(S_0)}\oplus Te(S_0)$ in $T(P\times\R^r)|_{e(S_0)}$ (with respect to some Riemannian metric). Choose a small neighbourhood $A$ of $e(S_0)$ in $\exp(\VV^\perp)$ (where $\exp$ denotes the exponential map of $P\times\R^r$) and for all $a\in A$ let $L_a$ and $L'_a$ be the intersections of a small neighbourhood of $a$ and the leaves of $\FF$ and $\FF'$ respectively.

If the neighbourhoods were chosen sufficiently small, then we are in the setting of \ref{tr3}, hence a diffeomorphism $\varphi\colon U\to U'$ arises (with the same notations as above). Note that $U$ and $U'$ are both neighbourhoods of $e(S_0)$, the map $\varphi$ fixes $e(S_0)$ and for all $a\in e(S_0)$ we have $d\varphi_a=\id_{T_a(P\times\R^m)}$. Observe that where the foliations $\FF$ and $\FF'$ initially coincide, this method just gives the identity for all $t\in[0,1]$. The isotopy we get this way can be extended to a diffeotopy of $P\times\R^r$ (by the isotopy extension theorem) and it takes the leaves of $\FF$ to the leaves of $\FF'$ around the image of $S_0$.

Next we repeat the same procedure around $e(S_1)$, the image of the next stratum, to get a new diffeotopy (that leaves a neighbourhood of $e(S_0)$ unchanged), and so on. In the end we obtain a diffeotopy of $P\times\R^r$ which turns $\FF$ into the vertical foliation $\FF'$ around the image of $M$ and does not change the vertical vector fields $\VV=\VV'$.
\end{prf}

\begin{thm}\label{taugxiemb}
For any set $\tau$ of stable singularities and any manifold $P^{n+k+m}$, if the number $r$ is sufficiently large (compared to $n$), then we have
$$\Cob_\tau^{G\oplus\xi}(P)\cong\Emb_\tau^{G\oplus\xi}(n,P\times\R^r).$$
\end{thm}

\begin{prf}
Take any number $r\ge2n+4$, so any manifold of dimension at most $n+1$ can be embedded into $\R^r$ uniquely up to isotopy. We will define two homomorphisms $\alpha\colon\Cob_\tau^{G\oplus\xi}(P)\to\Emb_\tau^{G\oplus\xi}(n,P\times\R^r)$ and $\beta\colon\Emb_\tau^{G\oplus\xi}(n,P\times\R^r)\to\Cob_\tau^{G\oplus\xi}(P)$ which will turn out to be each other's inverses.

\medskip\noindent I. \emph{Construction of $\alpha\colon\Cob_\tau^{G\oplus\xi}(P)\to\Emb_\tau^{G\oplus\xi}(n,P\times\R^r)$.}\medskip

For a $\tau^{G\oplus\xi}$-map $\tilde f\colon\xi_M\to P^{n+k}$ we can choose any embedding $i\colon M^n\hookrightarrow\R^r$ and define a $\tau^{G\oplus\xi}$-embedding $(e,\VV,\xi_M,\FF)$ as in example \ref{vert}. Define the map $\alpha$ to assign to the cobordism class of $f$ the cobordism class of $e=(e,\VV,\xi_M,\FF)$. In order to prove that $\alpha$ is well-defined, we have to show that the cobordism class of $e$ does not depend on the choice of the embedding $i$ and the representative of the class $[f]$.

\medskip\begin{sclaim}
If $i_0\colon M\hookrightarrow\R^r$ and $i_1\colon M\hookrightarrow\R^r$ are two embeddings and the above method assigns to them the $\tau^{G\oplus\xi}$-embeddings $e_0=(e_0,\VV_0,\xi_M,\mathscr{F}_0)$ and $e_1=(e_1,\VV_1,\xi_M,\mathscr{F}_1)$ respectively, then $e_0\sim e_1$.
\end{sclaim}

\begin{sprf}
Because of the dimension condition, $i_0$ and $i_1$ can be connected by an isotopy $i_t~(t\in[0,1])$. We define a $\tau^{G\oplus\xi}$-embedding
$$E\colon M\times[0,1]\hookrightarrow P\times\R^r\times[0,1];~(p,t)\mapsto (f(p),i_t(p),t)$$
(again with the horizontal $\xi_M$ and the vertical vector fields and foliation), which is precisely a cobordism between $e_0$ and $e_1$.
\end{sprf}

\begin{sclaim}
If $\tilde f_0\colon\xi_{M_0}\to P$ and $\tilde f_1\colon\xi_{M_1}\to P$ are cobordant $\tau^{G\oplus\xi}$-maps and the above method assigns to them the $\tau^{G\oplus\xi}$-embeddings $e_0=(e_0,\VV_0,\xi_{M_0},\mathscr{F}_0)$ and $e_1=(e_1,\VV_1,$ $\xi_{M_1},\mathscr{F}_1)$ respectively, then $e_0\sim e_1$.
\end{sclaim}

\begin{sprf}
Let $\tilde F\colon\xi_W\to P\times[0,1]$ be a cobordism between $\tilde f_0$ and $\tilde f_1$ (where $\xi_W$ is a vector bundle over a manifold $W^{n+1}$). Again by the dimension condition, the embedding $i_0\sqcup i_1\colon M_0\sqcup M_1=\partial W\hookrightarrow\R^r$ extends to an embedding $I\colon W\hookrightarrow\R^r$. Hence the map $E:=F\times I$ is a $\tau^{G\oplus\xi}$-embedding of $W$ into $P\times\R^r\times[0,1]$ (the vector fields and foliation are again vertical and $\xi_W$ is horizontal) and it is easy to see that this is a cobordism between $e_0$ and $e_1$.
\end{sprf}

\noindent II. \emph{Construction of $\beta\colon\Emb_\tau^{G\oplus\xi}(n,P\times\R^r)\to\Cob_\tau^{G\oplus\xi}(P)$.}\medskip

If $e=(e,\VV,\xi_M,\FF)$ is a $\tau^{G\oplus\xi}$-embedding of a manifold $M^n$ into $P^{n+k+m}\times\R^r$, then by the above lemma we obtain a diffeotopy of $P\times\R^r$ that turns $\VV$ and $\FF$ vertical. A diffeotopy of $P\times\R^r$ also yields a cobordism of $\tau^{G\oplus\xi}$-embeddings, hence we can assume that $\VV$ and $\FF$ were initially vertical. Now we can define $\beta$ to assign to the cobordism class of $e$ the cobordism class of the $\tau^{G\oplus\xi}$-map $\tilde f\colon\xi_M\to P$ for which $\tilde f|_M=\pr_P\circ e$ and $\tilde f|_{\xi_M|_p}=\pr_P\circ\exp_p$ for all $p\in e(M)$ (where $\pr_P$ denotes the projection to $P$ and $\exp_p$ now denotes the exponential of the leaf of $\FF$ at $p$). In order to prove that $\beta$ is well-defined, we have to show that the cobordism class of $\tilde f$ does not depend on the choice of the representative of the cobordism class $[e]$.

\medskip\begin{sclaim}
If $e_0=(e_0,\VV_0,\xi_{M_0},\mathscr{F}_0)$ and $e_1=(e_1,\VV_1,\xi_{M_1},\mathscr{F}_1)$ are cobordant $\tau^{G\oplus\xi}$-embed\-dings of the manifolds $M_0$ and $M_1$ respectively into $P\times\R^r$ and the above method assigns to them the $\tau^{G\oplus\xi}$-maps $\tilde f_0\colon\xi_{M_0}\to P$ and $\tilde f_1\colon\xi_{M_1}\to P$ respectively, then $f_0\sim f_1$.
\end{sclaim}

\begin{sprf}
We applied a diffeotopy $\varphi_t^i~(t\in[0,1])$ of $P\times\R^r\times\{i\}$ to turn the vector fields $\VV_i$ and foliation $\FF_i$ vertical (for $i=0,1$), this way we obtained the $\tau^{G\oplus\xi}$-map $\tilde f_i\colon\xi_{M_i}\to P\times\{i\}$. If $e_0$ and $e_1$ are connected by a cobordism $E=(E,\mathscr{U},\xi_W,\mathscr{G})$, which is a $\tau^{G\oplus\xi}$-embedding of a manifold $W^{n+1}$ into $P\times\R^r\times[0,1]$, then we can apply (the relative version of) the above lemma to obtain a diffeotopy $\Phi_t~(t\in[0,1])$ of $P\times\R^r\times[0,1]$ that extends the given diffeotopies $\varphi_t^0$ and $\varphi_t^1$ on the boundary and turns the vector fields $\mathscr{U}$ and the foliation $\mathscr{G}$ vertical. Now composing $E$ with the final diffeomorphism $\Phi_1$ and the projection to $P\times[0,1]$ as above, we obtain a $\tau^{G\oplus\xi}$-cobordism $\tilde F\colon\xi_W\to P\times[0,1]$ between $f_0$ and $f_1$ for which $\tilde F|_W=\pr_{P\times[0,1]}\circ\Phi_1\circ E$ and $\tilde F|_{\xi_W|_p}=\pr_{P\times[0,1]}\circ\Phi_1\circ\exp_p$ for all $p\in E(W)$.
\end{sprf}

 The constructions of $\alpha$ and $\beta$ imply also that they are homomorphisms and by design $\beta$ is the inverse of $\alpha$, hence they are both isomophisms between $\Cob_\tau^{G\oplus\xi}(P)$ and $\Emb_\tau^{G\oplus\xi}(n,P\times\R^r)$.
\end{prf}\medskip



Now we can describe the classifying space $X_\tau^{G\oplus\xi}$ for which we shall need the pullback of $\xi\to BG$ by the fibration $\pi\colon K_\tau^G\to BG$ (see definition \ref{kazsp}).

\begin{thm}\label{classtauxi}
If $\tau$ is a set of stable singularities, then we have $X_\tau^{G\oplus\xi}\cong\Gamma T(\nu_\tau^G\oplus\pi^*\xi)$.
\end{thm}

\begin{prf}
We need to prove that cobordisms of $\tau^{G\oplus\xi}$-maps to an arbitrary fixed manifold $P^{n+k+m}$ bijectively correspond to homotopy classes of maps $\cpt P\to\Gamma T(\nu_\tau^G\oplus\pi^*\xi)$. Since we have $\Cob_\tau^{G\oplus\xi}(P)\cong\Emb_\tau^{G\oplus\xi}(n,P\times\R^r)$ and $[\cpt P,\Gamma T(\nu_\tau^G\oplus\pi^*\xi)]\cong[S^r\cpt P,S^rT(\nu_\tau^G\oplus\pi^*\xi)]$ for sufficiently large numbers $r$, it is enough to prove
$$\Emb_\tau^{G\oplus\xi}(n,P\times\R^r)\cong[S^r\cpt P,S^rT(\nu_\tau^G\oplus\pi^*\xi)]$$
if $r$ is large enough. As in the proof of the previous theorem we will define two homomorphisms $\alpha\colon\Emb_\tau^{G\oplus\xi}(n,P\times\R^r)\to[S^r\cpt P,S^rT(\nu_\tau^G\oplus\pi^*\xi)]$ and $\beta\colon[S^r\cpt P,S^rT(\nu_\tau^G\oplus\pi^*\xi)]\to\Emb_\tau^{G\oplus\xi}(n,P\times\R^r)$ which are each other's inverses as we shall see.

\medskip\noindent I. \emph{Construction of $\alpha\colon\Emb_\tau^{G\oplus\xi}(n,P\times\R^r)\to[S^r\cpt P,S^rT(\nu_\tau^G\oplus\pi^*\xi)]$.}\medskip

Let $e=(e,\VV,\xi_M,\FF)$ be a $\tau^{G\oplus\xi}$-embedding of a manifold $M^n$ into $P\times\R^r$. The normal bundle of $e$ is of the form $\nu_e=\nu'\oplus\xi_M$ and by definition it is induced from $\gamma_{k+r}^G\oplus\xi|_{BG(k+r)}$ by a map $\mu_e\colon M\to BG(k+r)$. The embedding $e$ can be written as the composition
$$M\xra iD\xi_M\xra jD\nu_e\into P\times\R^r$$
where the disk bundles $D\xi_M$ and $D\nu_e$ are viewed as closed embedded tubular neighbourhoods with $j^{-1}(S\nu_e)=S\xi_M$ and the fibres of $D\xi_M$ are identified with the intersections of the leaves of $\FF$ with the tubular neighbourhood $D\xi_M$.

Now the vector fields $\VV$ can be extended to $D\xi_M$ by translation along the fibres and if $L$ is the leaf of $\FF$ at a point $p\in e(M)$, then we can use its exponential to define a foliation of dimension $r$ on $L$ which is tangent to the extended $\VV$ along the disk $D_p\xi_M$. Applying this for all points $p\in e(M)$ we obtain the structure of a $\tau^G$-embedding on $j$. Hence there are inducing maps $\kappa_j\colon D\xi_M\to K_\tau^G$ and $\tilde\kappa_j\colon D\nu_e\to S^rT\nu_\tau^G$ that pull back $j$ from the embedding $K_\tau^G\into S^rT\nu_\tau^G$ (again we are considering such a large number $r$ and such a finite dimensional approximation of the Kazarian space that the space $S^rT\nu_\tau^G$ exists).

The normal bundle of $e(M)$ in the disk bundle $D\xi_M$ is induced by the $\mu_e$ above and this all fits into a diagram
$$\xymatrix@R=1pc@C=1pc{
&&& D\nu_e\ar[rr]^{\tilde\kappa_j}\ar@{^(->}[dd] && S^rT\nu_\tau^G\ar[rr] && T\gamma_{k+r}^G \\
T\xi|_{BG(k+r)} && D\xi_M\ar[ll]\ar@{^(->}[ur]^j\ar[rr]^(.67){\kappa_j} && K_\tau^G\ar@{^(->}[ur]\ar[rr]^\pi && BG(k+r)\ar@{^(->}[ur] &\\
&&& P\times\R^r &&&&\\
BG(k+r)\ar@{^(->}[uu] && M\ar[ll]_(.33){\mu_e}\ar@{^(->}[uu]^i\ar@{^(->}[ur]^e &&&&&
}$$
Since $\mu_e$ induces the normal bundle of $e$, we have in the above diagram $\mu_e=\pi\circ\kappa_j\circ i$, that is, the inducing map $\mu_e$ factors through the Kazarian space $K_\tau^G$. But then we have a diagram
$$\xymatrix@R=1pc@C=1pc{
& \mu_e^*\xi=\xi_M\ar[dd]\ar[rr] && \pi^*\xi\ar[dd]\ar[rr] && \xi|_{BG(k+r)}\ar[dd] \\
\mu_e^*\gamma_{k+r}^G=\nu'\ar[dr]\ar[rr] && \pi^*\gamma_{k+r}^G=\nu_\tau^G\oplus\varepsilon^r\ar[dr]\ar[rr] && \gamma_{k+r}^G\ar[dr] &\\
& M\ar[rr]^{\kappa_e}\ar@/_1pc/[rrrr]_{\mu_e} && K_\tau^G\ar[rr]^(.4)\pi && BG(k+r) \\
}$$
where all squares are pullback squares.

Considering now the sum $\gamma_{k+r}\oplus\xi|_{BG(k+r)}$, we get $\nu_e=\kappa_e^*\pi^*(\gamma_{k+r}\oplus\xi|_{BG(k+r)})=\kappa_e^*(\nu_\tau^G\oplus\varepsilon^r\oplus\pi^*\xi)$, hence the embedding of $M$ into $P\times\R^r$ is induced by a diagram
$$\xymatrix{
S^r\cpt P\ar[r]^(.35){\tilde\kappa_e} & S^rT(\nu_\tau^G\oplus\pi^*\xi) \\
M\ar@{^(->}[u]^e\ar[r]^{\kappa_e} & K_\tau^G\ar@{^(->}[u]
}$$
This way we can assign to the $\tau^{G\oplus\xi}$-embedding $e$ the map $\tilde\kappa_e$ and since we can apply this construction to cobordisms as well, the assignment $\alpha[e]:=[\tilde\kappa_e]$ is well-defined.

\medskip\noindent II. \emph{Construction of $\beta\colon[S^r\cpt P,S^rT(\nu_\tau^G\oplus\pi^*\xi)]\to\Emb_\tau^{G\oplus\xi}(n,P\times\R^r)$.}\medskip

Suppose we have a map $\tilde\kappa\colon S^r\cpt P\to S^rT(\nu_\tau^G\oplus\pi^*\xi)$ and put $M^n:=\tilde\kappa^{-1}(K_\tau^G)$ and $\kappa:=\tilde\kappa|_M$. Now the normal bundle of $M$ in $P\times\R^r$ is of the form $\nu'\oplus\xi_M$ where $\nu'$ and $\xi_M$ are induced from $\gamma_{k+r}^G$ and $\xi|_{BG(k+r)}$ respectively by the map $\pi\circ\kappa\colon M\to BG(k+r)$.

If we take the preimage of the disk bundle $D\xi|_{BG(k+r)}\subset D(\gamma_{k+r}^G\oplus\xi|_{BG(k+r)})$ in $P\times\R^r$ (of a sufficiently small radius), then we get an embedding $j\colon D\xi_M\into D\nu_M$ where $\nu_M$ is the normal bundle of $M$ and the disk bundles are again considered as closed tubular neighbourhoods. Now we can project $D\nu_M$ to the orthogonal complement of $D\xi_M$ in it along the fibres of $D\xi_M$ and if we denote this projection by $\rho$, then the composition $\tilde\kappa\circ\rho$ is a map of $D\nu_M$ to $S^rT\nu_\tau^G$ such that the preimage of $K_\tau^G$ is $D\xi_M$, hence it gives the structure of a $\tau^G$-embedding to $j$.

We obtain vector fields $\VV=(v_1,\ldots,v_r)$ along $M$ by restricting the vector fields along $D\xi_M$. We also have a foliation $\GG$ of dimension $r$ on a neighbourhood of $M$ that is tangent to $\VV$ along $M$ and the projection $\rho$ maps each leaf of $\GG$ again into a leaf of $\GG$ because of the definition of the inducing map of the $\tau^G$-embedding $j$. So we can define a foliation $\FF$ of dimension $r+m$ by defining the leaf of $\FF$ at a point $p$ to be the preimage of the leaf of $\GG$ at $p$ under $\rho$.

This yields the structure of a $\tau^{G\oplus\xi}$-embedding $e=(e,\VV,\xi_M,\FF)$ on the embedding of $M$ into $P\times\R^r$. The same method assigns to a homotopy of $\tilde\kappa\colon S^r\cpt P\to S^rT(\nu_\tau^G\oplus\pi^*\xi)$ a cobordism of $e$, hence we can define $\beta[\tilde\kappa]:=[e]$.

\medskip
The above constructions imply that $\alpha$ and $\beta$ are homomorphisms and also that they are inverses of each other, hence we have proved $\Emb_\tau^{G\oplus\xi}(n,P\times\R^r)\cong[S^r\cpt P,S^rT(\nu_\tau^G\oplus\pi^*\xi)]$ and this is what we wanted.
\end{prf}

\section{Another long exact sequence}\label{another}

We shall apply the results of the section above in the case $G=\O$ and $\xi=2\Lambda\gamma=\Lambda\gamma\oplus\Lambda\gamma$ (recall that $\Lambda\gamma$ was the line bundle over $B\O$ which induces $\Lambda\gamma_n^\O$ over $B\O(n)$ for all $n$). The pullback of $\Lambda\gamma$ over the Kazarian space $K_\tau^\O$ is then $\Lambda\nu_\tau^\O$ 
and so the classifying space of $\tau^{\O\oplus2\Lambda\gamma}$-cobordisms is $\Gamma T(\nu_\tau^\O\oplus2\Lambda\nu_\tau^\O)$.

\begin{rmk}
Note that informally a $\tau^{\O\oplus2\Lambda\gamma}$-map can be thought of as a ``$\tau$-map'' $f\colon M^n\to P^{n+k+2}$ with normal bundle of the form $\nu_f=\nu'\oplus2\Lambda\nu'$. 
Now this $\Lambda\nu'$ coincides with $\Lambda\nu_f$ hence $f$ is a map with such a virtual normal bundle from which two (non-virtual) line bundles isomorphic to its determinant bundle split off.
\end{rmk}

The rest of part \ref{p2} is mainly devoted to the proof of theorem \ref{t2}. This will be quite similar to the proof of theorem \ref{t1} in part \ref{p1}.

\begin{thm}\label{ses}
For any set $\tau$ of stable singularities and any manifold $Q^q$ there is a long exact sequence
\begin{alignat*}2
\ldots&\xra{\psi'_{m+1}}\Cob_\tau^\i(Q\times\R^m)\xra{\varphi'_m}\Cob^\O_\tau(Q\times\R^m)\xra{\chi'_m}\Cob_\tau^{\O\oplus2\Lambda\gamma}(Q\times\R^m)\xra{\psi'_m} \\
&\xra{\psi'_m}\Cob_\tau^\i(Q\times\R^{m-1})\to\ldots
\end{alignat*}
\end{thm}

\begin{prf}
Recall the pullback diagram
$$\xymatrix{
B\i\ar[r]\ar@{^(->}[d] & S^1\ar@{^(->}[d]^i \\
B\O\ar[r]^{w_1} & \RP^\infty 
}$$
from definition \ref{bidef}. The complement of $S^1$ in $\RP^\infty$ deformation retracts to $\RP^{\infty-2}$ which is also a deformation retract of $\RP^\infty$, hence by pulling back its embedding by the fibration $w_1\colon B\O\to\RP^\infty$ to the embedding of a space $B$ we get a homotopy equivalence as shown on the diagram
$$\xymatrix{
B\O\ar[r]^{w_1} & \RP^\infty \\
B\ar[r]\ar@{_(->}[u]^\cong & \RP^{\infty-2}\ar@{_(->}[u]^\cong
}$$

Next we pull back the bundle $K_\tau^\O\to B\O$ by the homotopy equivalence $B\to B\O$ to a bundle $K\to B$, thus we have $K\cong K_\tau^\O$. Then the pullback of the universal virtual normal bundle $\nu_\tau^\O$ over $K$, which will be denoted $\nu$, is stably isomorphic to $\nu_\tau^\O$.

Take the virtual bundle $\nu_\tau^\i$ over $K_\tau^\i\subset K_\tau^\O$ (which is the restriction of $\nu_\tau^\O$) and the cofibration
$$S^nT\nu_\tau^\i(n)\into S^nT\nu_\tau^\O(n)\to S^nT\nu_\tau^\O(n)/S^nT\nu_\tau^\i(n)$$
for any $n$. The normal bundle of $K\subset K_\tau^\O$ is induced from the normal bundle of $\RP^{\infty-2}\subset\RP^{\infty}$ by the composition of $w_1\colon B\to\RP^{\infty-2}$ with the fibration $K\to B$. The normal bundle of $\RP^{\infty-2}$ is $2\gamma_1^\O$, this induces $2\Lambda\gamma$ over $B$ and this finally induces $2\Lambda\nu$ over the Kazarian space $K$, hence the cofibration above has the form
$$S^nT\nu_\tau^\i(n)\into S^nT\nu_\tau^\O(n)\to S^nT(\nu_\tau^\O(n)\oplus2\Lambda\nu_\tau^\O(n)).$$

Now applying the functor $\Omega^{n+m}\Gamma$ to the Puppe sequence of this cofibration we get a sequence of maps
$$\ldots\to\Omega^m\Gamma T\nu_\tau^\i(n)\to\Omega^m\Gamma T\nu_\tau^\O(n)\to\Omega^m\Gamma T(\nu_\tau^\O(n)\oplus2\Lambda\nu_\tau^\O(n))\to\Omega^{m-1}\Gamma T\nu_\tau^\i(n)\to\ldots$$
This sequence is infinite to the right by construction, but it is also infinite to the left since the number $n$ could be arbitrary and we get the same maps by applying $\Omega^{n+m}\Gamma$ to the $n$'th suspensions as by applying $\Omega^{n+m+1}\Gamma$ to the $(n+1)$'st suspensions. We can now converge with $n$ to infinity since the maps in this sequence commute with the natural maps induced by the inclusions $K_\tau^\sigma(n)\subset K_\tau^\sigma(n+1)$ (for $\sigma=\i,\O$) and so we get a sequence of maps
$$\ldots\to\Omega^m\Gamma T\nu_\tau^\i\to\Omega^m\Gamma T\nu_\tau^\O\to\Omega^m\Gamma T(\nu_\tau^\O\oplus2\Lambda\nu_\tau^\O)\to\Omega^{m-1}\Gamma T\nu_\tau^\i\to\ldots$$

If we then fix a manifold $Q$ and apply the functor $[\cpt Q,\cdot]$ to this sequence, then we obtain the long exact sequence of cobordism groups as claimed.
\end{prf}

\begin{rmk}
Similarly to the oriented case, if we had $Q=\R^q$, then the same long exact sequence could be obtained by turning the cofibration
$$T\nu_\tau^\i\into T\nu_\tau^\O\to T(\nu_\tau^\O\oplus2\Lambda\nu_\tau^\O)$$
(cf. remark \ref{cofib}) into the fibration
$$\Gamma T\nu_\tau^\O\xra{\Gamma T\nu_\tau^\i}\Gamma T(\nu_\tau^\O\oplus2\Lambda\nu_\tau^\O)$$
by the functor $\Gamma$ and taking its homotopy long exact sequence.
\end{rmk}

Observe that the restriction of $2\gamma_1^\O$ over the $1$-cell $S^1\subset\RP^\infty$ is trivial, hence the pullback of $2\Lambda\gamma\to B\O$ over $B\i$ will also be trivial which then induces the trivial bundle over the Kazarian space $K_\tau^\i$ as well. Thus analogously to the proof of theorem \ref{ses} we get a cofibration
$$S^2T\nu_\tau^\i\into T(\nu_\tau^\O\oplus2\Lambda\nu_\tau^\O)\to T(\nu_\tau^\O\oplus4\Lambda\nu_\tau^\O)$$
and iterating this process yields:

\begin{crly}
For any set $\tau$ of stable singularities and any manifold $Q^q$ and integer $r\ge0$ there is a long exact sequence
\begin{alignat*}2
\ldots&\to\Cob_\tau^\i(Q\times\R^{m-2r})\to\Cob^{\O\oplus2r\Lambda\gamma}_\tau(Q\times\R^m)\to\Cob_\tau^{\O\oplus2(r+1)\Lambda\gamma}(Q\times\R^m)\to \\
&\to\Cob_\tau^\i(Q\times\R^{m-2r-1})\to\ldots
\end{alignat*}
\end{crly}

\begin{rmk}
It would be tempting to try finding homomorphisms which complete the commutative diagram
$$\xymatrix@C=.75pc{
\ldots\ar[r] & \Cob_\tau^\i(Q\times\R^{m})\ar[d]^{\id}\ar[r] & \Cob^{\O\oplus2r\Lambda\gamma}_\tau(Q\times\R^{m+2r})\ar@{-->}[d]\ar[r] & \Cob_\tau^{\O\oplus2(r+1)\Lambda\gamma}(Q\times\R^{m+2r})\ar@{-->}[d]\ar[r] & \ldots \\
\ldots\ar[r] & \Cob_\tau^\i(Q\times\R^{m})\ar[r] & \Cob^{\O\oplus2s\Lambda\gamma}_\tau(Q\times\R^{m+2s})\ar[r] & \Cob_\tau^{\O\oplus2(s+1)\Lambda\gamma}(Q\times\R^{m+2s})\ar[r] & \ldots
}$$
with the dashed arrows for different numbers $r$ and $s$, however, it seems that such homomorphisms do not exist in general.
\end{rmk}

Now it remains from the proof of theorem \ref{t2} to describe the homomorphisms $\varphi'_m$ and $\chi'_m$ in the exact sequence in theorem \ref{ses}.

\section{Description of the homomorphisms $\varphi'_m$, $\chi'_m$ and $\psi'_m$}\label{phi'}

\begin{prop}
$\varphi'_m\colon\Cob_\tau^\i(Q\times\R^m)\to\Cob^\O_\tau(Q\times\R^m)$ is the forgetful homomorphism that assigns to the cobordism class of a Wall map $(f,w)$ the unoriented cobordism class of $f$.
\end{prop}

\begin{prf}
This follows immediately since the map between the classifying spaces is just the inclusion of $\Gamma T\nu_\tau^\i$ into $\Gamma T\nu_\tau^\O$ (see the proof of theorem \ref{ses}).
\end{prf}

\begin{prop}
$\chi'_m\colon\Cob_\tau^\O(Q\times\R^m)\to\Cob_\tau^{\O\oplus2\Lambda\gamma}(Q\times\R^m)$ assigns to a cobordism class $[f]$ the cobordism class of $f$ restricted to the Poincaré dual of $w_1(\nu_f)^2$ which can be represented by a $2$-codimensional submanifold uniquely up to cobordism, the restriction $f|_{\PD(w_1(\nu_f)^2)}$ has a normal $\O\oplus2\Lambda\gamma$-structure and its cobordism class only depends on the class of $f$, hence this assignment is well-defined.
\end{prop}

\begin{prf}
This is completely analogous to proposition \ref{pd}. For an unoriented cobordism class $[f]\in\Cob_\tau^\O(Q\times\R^m)$ represented by a $\tau$-map $f\colon M\to Q\times\R^m$ the Poincaré dual of
$$w_1(\nu_f)\colon M\to\RP^\infty$$
is represented by the preimage $w_1(\nu_f)^{-1}(\RP^{\infty-1})$. Then choosing two distinct hyperplanes $\RP^{\infty-1}_1$ and $\RP^{\infty-1}_2$ in $\RP^\infty$ the Poincaré dual of $w_1(\nu_f)^2$ is represented by
\begin{alignat*}2
w_1(\nu_f)^{-1}(\RP^{\infty-1}_1)\cap w_1(\nu_f)^{-1}(\RP^{\infty-1}_2)&=w_1(\nu_f)^{-1}(\RP^{\infty-1}_1\cap\RP^{\infty-1}_2)= \\
&=w_1(\nu_f)^{-1}(\RP^{\infty-2}).
\end{alignat*}
Now the same argument as in the proof of proposition \ref{pd} yields that $\chi'_m$ assigns to the cobordism class of $f$ the class of $f|_{w_1(\nu_f)^{-1}(\RP^{\infty-2})}=f|_{\PD(w_1(\nu_f)^2)}$.
\end{prf}\medskip

This finishes the proof of theorem \ref{t2}.

\begin{rmk}\label{psi'}
Although there are no apparent interesting algebraic properties of the homomorphism $\psi'_m\colon\Cob_\tau^{\O\oplus2\Lambda\gamma}(Q\times\R^m)\to\Cob_\tau^\i(Q\times\R^{m-1})$, it has a nice geometric description. Let $q$ be the dimension of $Q$ and put $n:=q+m-k$. Observe that we have
\begin{alignat*}2
\Cob_\tau^{\O\oplus2\Lambda\gamma}(Q^q\times\R^m)&\cong\Emb_\tau^{\O\oplus2\Lambda\gamma}(n-2,Q\times\R^{m+r})\cong\\
&\cong\Emb^{\nu_\tau^\O\oplus\varepsilon^r\oplus2\Lambda\nu_\tau^\O}(n-2,Q\times\R^{m+r})
\end{alignat*}
and
\begin{alignat*}2
\Cob_\tau^\i(Q^q\times\R^{m-1})&\cong\Emb_\tau^{\i\oplus1}(n-1,Q\times\R^{m+r})\cong\\
&\cong\Emb^{\nu_\tau^\i\oplus\varepsilon^{r+1}}(n-1,Q\times\R^{m+r})
\end{alignat*}
where $r$ is sufficiently large, moreover, $\psi'_m$ is obtained as the boundary homomorphism in the Puppe sequence of classifying spaces. Hence we can apply lemma \ref{iotal} by setting the $K$, $A$, $\xi$, $B$ and $\zeta$ in the lemma to be $K_\tau^\O$, $K$, $2\Lambda\nu$, $K_\tau^\i$ and $\nu_\tau^\O\oplus\varepsilon^r$ respectively. This yields that for any cobordism class $[\tilde f]\in\Cob_\tau^{\O\oplus2\Lambda\gamma}(Q^q\times\R^m)$ its image $\psi'_m[\tilde f]$ is represented by the mapping of $S\xi_M$ (that is, the circle bundle of the $2\Lambda\gamma$-part of the normal bundle of $\tilde f|_M$) by the restriction of a representative of $\tilde f$, together with its natural outward normal vector field.

\end{rmk}

\begin{rmk}\label{gen2}
As in remark \ref{gen1} consider again the case $Q\times\R^m=\R^{n+k}$ where the codimension $k$ is large (compared to $n$). Then for any (non-empty) singularity set $\tau$ we have $\Cob_\tau^\i(\R^{n+k})=\WW_n$ and $\Cob_\tau^\O(\R^{n+k})=\NN_n$, moreover, we also have $\Cob_\tau^{\O\oplus2\Lambda\gamma}(\R^{n+k})=:\Cob_\tau^{\O\oplus2\Lambda\gamma}(\R^{n+k})=\NN_{n-2}$ since for any embedding $i\colon M^{n-2}\into\R^{n+k}$ with normal bundle $\nu_i=\nu'\oplus2\Lambda\nu'$ this $\Lambda\nu'$ coincides with $\Lambda\nu_i$ which is isomorphic to $\Lambda TM$, so this normal structure only depends on $M$ and if $k$ is large enough, $2\Lambda TM$ can be embedded into $\R^{n+k}$ uniquely up to isotopy.

Thus theorem \ref{t2} gives an exact sequence for which the portion where $n$ is sufficiently small (compared to $k$) looks like
$$\ldots\to\WW_n\to\NN_n\to\NN_{n-2}\to\WW_{n-1}\to\ldots$$
and if we increase $k$ with a fixed $n$, then the homomorphisms do not change which means that this sequence is infinite both to the right and to the left. Then the propositions above show that this sequence can be identified with that in \cite[theorem 4.3]{atiyah}. Since the codimension $k$ was assumed to be large enough, the manifold constructed in \cite[theorem 4.4]{atiyah} can also be mapped to $\R^{n+k}$ uniquely up to isotopy which now gives a splitting $\NN_{n-2}\to\NN_n$ of this long exact sequence yielding the classical exact sequence (\ref{ces2}), hence theorem \ref{t2} really generalises the sequence (\ref{ces2}). Later, in proposition \ref{n0} we shall show that such a splitting cannot exist in general.
\end{rmk}

\part{Consequences and applications}\label{p3}


In the following two sections we shall apply our two exact sequences for the simplest types of singularity sets, i.e. for $\tau=\{\Sigma^0\}$ (the case of immersions) and $\tau=\{\Sigma^0,\Sigma^{1,0}\}$, $\tau=\{\Sigma^0,\Sigma^{1,0},\Sigma^{1,1,0}\}$ and so on (the case of Morin maps). Our first goal will always be to describe how the endomorphism $\varphi_m$ in theorem \ref{t1} acts rationally. By \cite[proposition 90]{hosszu} we have
$$\Cob_\tau^\SO(P^{n+k})\otimes\Q\cong\displaystyle\bigoplus_{i=1}^{n+k}H_i(P;\Q)\otimes H^{i-k}(K_\tau^\SO;\Q),$$
in particular $\Cob_\tau^\SO(\R^{n+k})\otimes\Q\cong H^n(K_\tau^\SO;\Q)$, hence in order to understand $\varphi_m$ rationally we only have to know how the involution which induces $\iota$ acts on the Kazarian space $K_\tau^\SO$.

After this rational description we shall compute in both cases some cobordism groups of Wall $\tau$-maps for the above singularity sets $\tau$ which is interesting because Wall cobordism groups are what connect oriented and unoriented cobordism groups. When viewing manifolds abstractly, the classical sequences (\ref{ces1}) and (\ref{ces2}) yield an obstruction for an unoriented cobordism class can be representable by an orientable manifold, namely that it should be a Wall cobordism class. Now the forgetful homomorphism $\varphi'_m$ in theorem \ref{t2} is not always injective but otherwise we have the same property, namely for an unoriented $\tau$-cobordism class to be representable by an orientable $\tau$-map it should be the $\varphi'_m$-image of a Wall $\tau$-cobordism class.

\section{Immersions}

In this section we investigate the case $\tau=\{\Sigma^0\}$, i.e. $\tau$-maps are the $k$-codimensional immersions. We shall use the notation
$$\Imm^\sigma(n,k):=\Cob_{\{\Sigma^0\}}^\sigma(\R^{n+k})\quad\text{and}\quad\Imm^\sigma(n,P^{n+k}):=\Cob_{\{\Sigma^0\}}^\sigma(P^{n+k})$$
for any stable normal normal structure $\sigma$ and any manifold $P$ (except that for $\sigma=G\oplus\xi$ we decrease $n$ and increase $k$ by the rank of $\xi$). Now the oriented Kazarian space is $K_{\{\Sigma^0\}}^\SO=B\SO(k)$ and the universal normal bundle over it is $\gamma_k^\SO$ and we immediately obtain:

\begin{prop}\label{immq}
The endomorphism $\varphi_{n+k}\otimes\Q$ of $\Imm^\SO(n,k)\otimes\Q$ is the following:
\begin{enumerate}
\item if $k=2m$ is even and $q$ is the number of non-negative integers $a_0,a_1,\ldots,a_m$ such that $a_0$ is odd and $n=a_0k+\underset{i=1}{\overset m\sum}a_i4i$, then $\varphi_{n+k}\otimes\Q$ is trivial on $q$ generators and the multiplication by $2$ on the rest of the generators (in an appropriate basis),
\item if $k$ is odd, then $\varphi_{n+k}\otimes\Q$ is the multiplication by $2$.
\end{enumerate}
\end{prop}

\begin{prf}
We have
\begin{alignat*}2
\Imm^\SO(*,k)\otimes\Q&\cong H^*(B\SO(k);\Q)=\\
&=\begin{cases}
\Q[p_1,\ldots,p_m,e]/(e^2-p_m),&\text{if }k=2m\text{ is even}\\
\Q[p_1,\ldots,p_m],&\text{if }k=2m+1\text{ is odd}
\end{cases}
\end{alignat*}
where the $p_i\in H^{4i}(B\SO(k);\Q)$ are the Pontryagin classes and $e\in H^k(B\SO;\Q)$ is the Euler class of $\gamma_k^\SO$. The identification of $\Imm^\SO(n,k)\otimes\Q\cong\pi^s_{n+k+1}(ST\gamma_k^\SO)$ with the degree-$n$ part of this graded ring follows from the stable Hurewicz homomorphism (which is rationally iso), the universal coefficient theorem and the Thom isomorphism corresponding to the bundle $\gamma_k^\SO\oplus\varepsilon^1\to B\SO(k)$.

By proposition \ref{2} (see also corollary \ref{iotac}) the involution inducing $\iota$ acts on this vector bundle so that it inverts the summand $\varepsilon^1$ and reverses orientation on the summand $\gamma_k^\SO$, thus the Pontryagin classes $p_i$ are unchanged by its action and the Euler class $e$ is mapped to $-e$. Our claim easily follows from this.
\end{prf}\medskip

We can say more about these cobordism groups in those cases where the codimension $k$ is either small or relatively large.

\subsection*{Small codimensional immersions}

First, if the codimension is $1$, then the normal bundle of any immersion $f\colon M^n\imto P^{n+1}$ is induced by $w_1(\nu_f)\colon M^n\to\RP^\infty$ from $\gamma_1^\O$. The classifying spaces of $\Imm^\SO$, $\Imm^\i$, $\Imm^\O$ and $\Imm^{\O\oplus2\Lambda\gamma}$ are $\Gamma T\varepsilon^1$, $\Gamma T\gamma^\O_1|_{\RP^1}$, $\Gamma T\gamma^\O_1$ and $\Gamma T(3\gamma^\O_1)$ respectively and since we have $T(k\gamma_1^\O|_{\RP^m})=\RP^{m+k}/\RP^{k-1}$ (for all $k$ and $m$) these classifying spaces are $\Gamma S^1$, $\Gamma\RP^2$, $\Gamma\RP^\infty$ and $\Gamma(\RP^\infty/\RP^2)$ respectively.

\begin{rmk}
The group $\Imm^{\O\oplus2\Lambda\gamma}(n-2,P^{n+1})$ is the cobordism group of those immersions $f\colon M^{n-2}\imto P^{n+1}$ whose normal bundle splits to the sum of three identical line bundles, i.e. the group $\Imm^{\text{sfr}}(n-2,P)$ of $3$-codimensional skew-framed immersions which arise naturally in the study of framed immersions; see \cite{ae}.
\end{rmk}

Now if we view immersions to Eucledian spaces, then theorems \ref{t1} and \ref{t2} yield the long exact sequences
$$\ldots\to\pi^s(n)\to\pi^s(n)\to\pi^s_{n+1}(\RP^2)\to\pi^s(n-1)\to\ldots$$
and
$$\ldots\to\pi^s_{n+1}(\RP^2)\to\pi^s_{n+1}(\RP^\infty)\to\pi^s_{n+1}(\RP^\infty/\RP^2)\to\pi^s_n(\RP^2)\to\ldots$$
of stable homotopy groups (which, of course, could be both obtained in much easier ways too). This means that in this case our main theorems are not new, however it also shows the following important property of our two exact sequences:

\begin{prop}\label{n0}
The composition $\psi_m\circ\psi'_{m+1}$ of the homomorphisms $\psi'_{m+1}$ in theorem \ref{t2} and $\psi_m$ in theorem \ref{t1} is not always zero.
\end{prop}

This is interesting since the classical analogue of $\psi'_{m+1}$ is always zero yielding that the long exact sequence splits to short exact sequences of the form (\ref{ces2}). However, this proposition shows that the general sequence in theorem \ref{t2} does not split to short exact sequences.

\medskip\begin{prf}
Consider the case of $1$-codimensional immersions of $4$-manifolds. Then the combination of the above two exact sequences gives a diagram
$$\xymatrix{
& \pi^s_6(\RP^\infty/\RP^2)\ar[d]^{\psi'_6} &\\
\pi^s(4)\ar[r] & \pi^s_5(\RP^2)\ar[r]^{\psi_5}\ar[d] & \pi^s(3) \\
& \pi^s_5(\RP^\infty) &
}$$
with exact row and column. Since $\pi^s(4)$ and $\pi^s_5(\RP^\infty)$ are trivial (see \cite{liu}) we get that $\psi'_6$ is epi and $\psi_5$ is mono and since $\pi^s_5(\RP^2)$ is $\Z_2$ (see \cite{wu}) this means that $\psi_5\circ\psi'_6$ is non-zero.
\end{prf}\medskip

Now let us consider immersions of codimension $2$. In this case the classifying space of $\Imm^\SO$ is $\Gamma T\gamma_2^\SO$ and since $\gamma_2^\SO$ coincides with the tautological complex line bundle over $\CP^\infty$ this space is $\Gamma\CP^\infty$. The first few stable homotopy groups of $\CP^\infty$ were computed by Liulevicius \cite{liu} and Mosher \cite{mosh} and are as follows:
\begin{table}[H]\begin{center}\begin{tabular}{c||c|c|c|c|c|c|c|c|c|c|c|c}
$m$ & $1$ & $2$ & $3$ & $4$ & $5$ & $6$ & $7$ & $8$ & $9$ & $10$ & $11$ & $12$\\
\hline
$\pi^s_m(\CP^\infty)$ & $0$ & $\Z$ & $0$ & $\Z$ & $\Z_2$ & $\Z$ & $\Z_2$ & $\Z\oplus\Z_2$ & $\Z_8\oplus\Z_3$ & $\Z$ & $\Z_4$ & $\Z\oplus\Z_3$
\end{tabular}\end{center}\end{table}\vspace{-.5cm}\noindent
The action of the involution $\iota$ in theorem \ref{t1} for $\Imm^\SO(n,2)\cong\pi^s_{n+2}(\CP^\infty)$ now immediately follows for $n\le6$ from proposition \ref{immq}, hence we have:

\begin{prop}
The endomorphism $\varphi_{n+2}$ of $\Imm^\SO(n,2)$ is
\begin{enumerate}
\item $0$ on the free part for all $n\equiv2~(4)$ and $2\id$ on the free part for all $n\equiv0~(4)$,
\item $0$ on the torsion part for $0\le n\le6$.
\end{enumerate}
\end{prop} 

\begin{crly}
The cobordism groups $\Imm^\i(n,2)$ for $n\le6$ are
\begin{table}[H]\begin{center}\begin{tabular}{c||c|c|c|c|c|c|c}
$n$ & $0$ & $1$ & $2$ & $3$ & $4$ & $5$ & $6$ \\
\hline
$\Imm^\i(n,2)$ & $\Z_2$ & $0$ & $\Z$ & $\Z\oplus\Z_2$ & $\Z_2~?~\Z_2$ & $\Z_2$ & $\Z\oplus(\Z_2~?~\Z_2)$
\end{tabular}\end{center}\end{table}\vspace{-.5cm}\noindent
where $G~?~H$ denotes the existence of a short exact sequence $0\to G\to G~?~H\to H\to0$.
\end{crly} 

\begin{prf}
Use theorem \ref{t1} to obtain the exact sequence
$$0\to\coker\varphi_{n+2}\to\Imm^\i(n,2)\to\ker\varphi_{n+1}\to0$$
and apply the proposition above.
\end{prf}

\subsection*{Large codimensional immersions}

Let us now turn to cobordisms of $k$-codimensional immersions such that the dimension of the source manifolds is not much greater than $k$. The reason for this is that for $n$ close to $k$ the group $\Imm^\sigma(n,k)$ is close to the abstract cobordism group of $n$-manifolds with normal $\sigma$-structures; see the works of Koschorke \cite{kos}, Olk \cite{olk}, Pastor \cite{pas} and Li \cite{li} where these groups were computed for $k\le n\le k+2$ and $\sigma=\O,\SO$.

\begin{rmk}\label{wueiaieviaqhevahhuev}
For $n<k$ we have $\Imm^\SO(n,k)=\Omega_n$ and $\Imm^\O(n,k)=\NN_n$ since a cobordism between two $n$-manifolds is at most a $k$-manifold which, when generically mapped with codimension $k$, is immersed.
\end{rmk}

We shall use the natural forgetful homomorphisms
\begin{gather*}
\alpha^\SO_\imm\colon\Imm^\SO(n,k)\to\Omega_n,\quad\alpha^\i_\imm\colon\Imm^\i(n,k)\to\WW_n,\\
\alpha^\O_\imm\colon\Imm^\O(n,k)\to\NN_n\quad\text{and}\quad\alpha^{\O\oplus2\Lambda\gamma}_\imm\colon\Imm^{\O\oplus2\Lambda\gamma}(n-2,k+2)\to\NN_{n-2}
\end{gather*}
that assign to the cobordism class of an immersion $f\colon M^n\imto\R^{n+k}$ (or a germ $\tilde f\colon\xi_M\to\R^{n+k}$) with the appropriate normal structure the abstract cobordism class of $M$. These homomorphisms commute with the exact sequences in theorems \ref{t1} and \ref{t2} and the classical exact sequences (\ref{ces1}) and (\ref{ces2}), that is, we have commutative diagrams
$$\xymatrix@C=.75pc{
\ldots\ar[r] & \Imm^\SO(n,k)\ar[d]^{\alpha^\SO_\imm}\ar[r] & \Imm^\SO(n,k)\ar[d]^{\alpha^\SO_\imm}\ar[r] & \Imm^\i(n,k)\ar[d]^{\alpha^\i_\imm}\ar[r] & \Imm^\SO(n-1,k)\ar[d]^{\alpha^\SO_\imm}\ar[r] & \ldots \\
\ldots\ar[r] & \Omega_n\ar[r]^{2\id} & \Omega_n\ar[r] & \WW_n\ar[r] & \Omega_{n-1}\ar[r]^{2\id} & \ldots
}$$
and
$$\xymatrix@C=.75pc{
\ldots\ar[r] & \Imm^\i(n,k)\ar[d]^{\alpha^\i_\imm}\ar[r] & \Imm^\O(n,k)\ar[d]^{\alpha^\O_\imm}\ar[r] & \Imm^{\O\oplus2\Lambda\gamma}(n-2,k+2)\ar[d]^{\alpha^{\O\oplus2\Lambda\gamma}_\imm}\ar[r] & \Imm^\i(n-1,k)\ar[d]^{\alpha^\i_\imm}\ar[r] & \ldots \\
\ldots\ar[r]^0 & \WW_n\ar[r] & \NN_n\ar[r] & \NN_{n-2}\ar[r]^0 & \WW_{n-1}\ar[r] & \ldots
}$$
where the rows are exact.

\begin{lemma}\label{immo+}
The map $\alpha^{\O\oplus2\Lambda\gamma}_\imm$ is an isomorphism $\Imm^{\O\oplus2\Lambda\gamma}(n-2,k+2)\cong\NN_{n-2}$ for all $n\le k+1$.
\end{lemma}

\begin{prf}
We define the inverse of $\alpha^{\O\oplus2\Lambda\gamma}_\imm$ in the following way: for a cobordism class $[M]\in\NN_{n-2}$ represented by $M^{n-2}$ take the bundle $\xi_M:=2\Lambda TM$ and map it generically to $\R^{n+k}$; let the germ of this map along $M$ be $\tilde f$ and assign the cobordism class $[\tilde f]$ to $[M]$. Such a correspondence would be the inverse of $\alpha^{\O\oplus2\Lambda\gamma}_\imm$ if it was well-defined since $\Lambda TM$ is also the determinant bundle of the stabilisation of $\nu_{(\tilde f|_M)}$. So we only have to prove that $\tilde f$ is an immersion germ and its cobordism class only depends on that of $M$.

Let $N^{n-2}$ be another representative of $[M]$ and let $W^{n-1}$ be a compact manifold with boundary such that $\partial W=M\sqcup N$. Then the restriction of $\xi_W=2\Lambda TW$ to the boundary gives $\xi_M=2\Lambda TM$ and $\xi_N=2\Lambda TN$ and if $\tilde f$ and $\tilde g$ are their generic map germs to $\R^{n+k}$ as defined above, then we can extend their representatives $f\colon\xi_M\to\R^{n+k}\times\{0\}$ and $g\colon\xi_N\to\R^{n+k}\times\{1\}$ by a generic map
$$F\colon\xi_W\to\R^{n+k}\times[0,1].$$

By the dimension condition such a map is stable. Its singular locus $\Sigma(F):=\{p\in W\mid\rk dF_p<n-1\}$ is a submanifold of $\xi_W$ and a straightforward computation with jet bundles shows that we have $\codim\Sigma(F)=k+1$. 
Now since the dimension of $W$ is at most $k$ we generically have $W\cap\Sigma(F)=\varnothing$, that is, the map $F$ in a small neighbourhood of the zero-section is an immersion. This means that the germ $\tilde F$ of $F$ along $W$ satisfies both conditions in definition \ref{deftauxi}, hence both $\tilde f$ and $\tilde g$ are $\O\oplus2\Lambda\gamma$-immersions and $\tilde F$ is an immersed $\O\oplus2\Lambda\gamma$-cobordism between them. This finishes the proof.
\end{prf}

\begin{prop}
For all $k\ge1$ we have
\begin{enumerate}
\item $\Imm^\i(k,k)\cong\WW_k\oplus\Z$ if $k$ is even,
\item $\Imm^\i(k,k)\cong\WW_k\oplus\Z_2$ if $k$ is odd.
\end{enumerate}
\end{prop}

\begin{prf}
By Koschorke \cite[theorem 10.8]{kos} and Pastor \cite[theorem 3.1]{pas} we have
$$\Imm^\O(k,k)=\NN_k\oplus G\quad\text{and}\quad\Imm^\SO(k,k)=\Omega_k\oplus G$$
where $G$ denotes $\Z$ if $k$ is even and $\Z_2$ if $k$ is odd. Now by proposition \ref{immq} and lemma \ref{immo+} the two diagrams above lemma \ref{immo+} take the form
$$\xymatrix{
0\ar[r] & \frac{\Omega_k}{2\Omega_k}\oplus G\ar[d]^{\id\oplus0}\ar[r] & \Imm^\i(k,k)\ar[d]^{\alpha^\i_\imm}\ar[r] & T_2(\Omega_{k-1})\ar[d]^{\id}\ar[r] & 0 \\
0\ar[r] & \frac{\Omega_k}{2\Omega_k}\ar[r] & \WW_k\ar[r] & T_2(\Omega_{k-1})\ar[r] & 0
}$$
and
$$\xymatrix{
& \Imm^\i(n,k)\ar[d]^{\alpha^\i_\imm}\ar[r] & \NN_k\oplus G\ar[d]^{\id\oplus0}\ar[r] & \NN_{k-2}\ar[d]^{\id}\ar[r] & 0 \\
0\ar[r] & \WW_k\ar[r] & \NN_k\ar[r] & \NN_{k-2}\ar[r] & 0
}$$
where $T_2(\Omega_{k-1})$ is the $2$-torsion subgroup of $\Omega_{k-1}$. The $5$-lemma (or the snake lemma) applied to the first of these diagrams implies that $\alpha^\i_\imm$ is epi and its kernel is $G$ and the second diagram implies that this kernel is also a direct summand.
\end{prf}

\begin{rmk}
Actually the proposition above could also be proved without using theorems \ref{t1} and \ref{t2} just by noting that the kernels of $\alpha^\SO_\imm$ and $\alpha^\O_\imm$ are generated by the same object: the double point set (up to cobordism) of an immersion of $S^k$ with one double point.
\end{rmk}

\begin{prop}
For all $k\ge1$ we have
\begin{enumerate}
\item\label{iw1} $\Imm^\i(k+1,k)\cong\WW_{k+1}\oplus\Z_4$ if $k\equiv1~(4)$,
\item\label{iw2} $\Imm^\i(k+1,k)\cong\WW_{k+1}\oplus\Z\oplus\Z_2$ if $k\equiv2~(4)$,
\item\label{iw3} $\Imm^\i(k+1,k)\cong\WW_{k+1}\oplus\Z_2$ if $k\equiv3~(4)$,
\item\label{iw4} $\Imm^\i(k+1,k)\cong\WW_{k+1}\oplus\Z\oplus\Z_2\oplus\Z_2$ if $k\equiv0~(4)$.
\end{enumerate}
\end{prop}

\begin{prf}
The first thing to note (similarly to the proof above) is that by lemma \ref{immo+} and the snake lemma applied to the second diagram above it we have an epimorphism $\ker\alpha^\i_\imm\onto\ker\alpha^\O_\imm$. Moreover, if $k+1$ is not a power of $2$ and $k\ne1$, then by \cite[theorem 10.8]{kos} $\ker\alpha^\O_\imm$ is a direct summand in $\Imm^\O(k+1,k)$ and so if the image of $[f]\in\ker\alpha^\i_\imm$ is non-zero in $\ker\alpha^\O_\imm$, then $[f]$ is independent in $\Imm^\i(k+1,k)$ of the subgroup generated by the $\alpha^\i_\imm$-preimage of $\WW_{k+1}\setminus\{0\}$.

In the following we shall always use the first diagram above lemma \ref{immo+} which yields (again by the $5$-lemma) that $\alpha^\i_\imm$ is epi in all cases. We will also use that by proposition \ref{immq} the endomorphism $\varphi_{2k}$ is $0$ on the direct complement of $\Omega_k$ in $\Imm^\SO(k,k)$ and we obtain the forms of $\Imm^\SO(k+1,k)$ and $\Imm^\O(k+1,k)$ from \cite[theorem 6]{li} and \cite[theorem 10.8]{kos} respectively.

\medskip\noindent\textit{Proof of \ref{iw1}.\enspace\ignorespaces}
If $k\equiv1~(4)$, we have
$$\xymatrix{
0\ar[r] & \Omega_{k+1}\oplus\Z_2\ar[d]^{\id\oplus0}\ar[r] & \Imm^\i(k+1,k)\ar[d]^{\alpha^\i_\imm}\ar[r] & \Omega_k\oplus\Z_2\ar[d]^{\id\oplus0}\ar[r] & 0 \\
0\ar[r] & \Omega_{k+1}\ar[r] & \WW_{k+1}\ar[r] & \Omega_k\ar[r] & 0
}$$
hence $\ker\alpha^\i_\imm$ is either $\Z_4$ or $\Z_2\oplus\Z_2$. But $\ker\alpha^\i_\imm$ maps onto $\ker\alpha^\O_\imm\cong\Z_4$, so it can only be $\Z_4$ which can only be a direct summand.

\medskip\noindent\textit{Proof of \ref{iw2}.\enspace\ignorespaces}
If $k\equiv2~(4)$ and $k+2$ is not a power of $2$, then we have $\Imm^\SO(k+1,k)\cong\Omega_{k+1}\oplus\Z_4$ and by \cite[p. 472]{li} the involution $\iota$ in theorem \ref{t1} is the identity, hence the cokernel of $\varphi_{2k+1}$ is $\Omega_{k+1}\oplus\Z_2$. If $k+2$ is a power of $2$, then we have $\Imm^\SO(k+1,k)\cong\Omega_{k+1}\oplus\Z_2$ and so in both cases we get the diagram
$$\xymatrix{
0\ar[r] & \Omega_{k+1}\oplus\Z_2\ar[d]^{\id\oplus0}\ar[r] & \Imm^\i(k+1,k)\ar[d]^{\alpha^\i_\imm}\ar[r] & \Omega_k\oplus\Z\ar[d]^{\id\oplus0}\ar[r] & 0 \\
0\ar[r] & \Omega_{k+1}\ar[r] & \WW_{k+1}\ar[r] & \Omega_k\ar[r] & 0
}$$
which implies that we have $\ker\alpha^\i_\imm\cong\Z\oplus\Z_2$. The $\Z$ part cannot be anything else than a direct summand in $\Imm^\i(k+1,k)$ but the $\Z_2$ part is also a direct summand since by the proof of \cite[theorem 3.2]{pas} it is generated by a cobordism class $[f]$ which maps to the generator of the $\Z_2$ part of $\Imm^\SO(k+1,k+1)\cong\Omega_{k+1}\oplus\Z_2$ under the homomorphism induced by the inclusion $\R^{2k+1}\subset\R^{2k+2}$ and this is non-zero in $\Imm^\O(k+1,k+1)$ as well, hence the class of $f$ in $\Imm^\O(k+1,k)$ cannot vanish either.

\medskip\noindent\textit{Proof of \ref{iw3}.\enspace\ignorespaces}
If $k\equiv3~(4)$, we have
$$\xymatrix{
0\ar[r] & \frac{\Omega_{k+1}}{2\Omega_{k+1}}\ar[d]^{\id}\ar[r] & \Imm^\i(k+1,k)\ar[d]^{\alpha^\i_\imm}\ar[r] & \Omega_k\oplus\Z_2\ar[d]^{\id\oplus0}\ar[r] & 0 \\
0\ar[r] & \frac{\Omega_{k+1}}{2\Omega_{k+1}}\ar[r] & \WW_{k+1}\ar[r] & \Omega_k\ar[r] & 0
}$$
hence $\ker\alpha^\i_\imm$ is $\Z_2$ which can only be a direct summand (even if $k+1$ is a power of $2$).

\medskip\noindent\textit{Proof of \ref{iw4}.\enspace\ignorespaces}
If $k\equiv0~(4)$, we have
$$\xymatrix{
\Omega_{k+1}\oplus\Z_2\oplus\Z_2\ar[r]^{\id+\iota} & \Omega_{k+1}\oplus\Z_2\oplus\Z_2\ar[d]^{\id\oplus0\oplus0}\ar[r] & \Imm^\i(k+1,k)\ar[d]^{\alpha^\i_\imm}\ar[r] & T_2(\Omega_k)\oplus\Z\ar[d]^{\id\oplus0}\ar[r] & 0 \\
0\ar[r] & \Omega_{k+1}\ar[r] & \WW_{k+1}\ar[r] & T_2(\Omega_k)\ar[r] & 0
}$$
and so $\ker\alpha^\i_\imm$ is either $\Z\oplus\Z_2$ or $\Z\oplus\Z_2\oplus\Z_2$ depending on whether the $\iota$ in theorem \ref{t1} swaps the two $\Z_2$ summands in $\Omega_{k+1}\oplus\Z_2\oplus\Z_2$ or not. Now again the proof of \cite[theorem 3.2]{pas} shows that one of these summands $\Z_2$ is generated by a cobordism class $[f]$ which maps to the generator of the $\Z_2$ part of $\Imm^\SO(k+1,k+1)\cong\Omega_{k+1}\oplus\Z_2$ under the homomorphism induced by the inclusion $\R^{2k+1}\subset\R^{2k+2}$. Since the reflection to a hyperplane commutes with this inclusion we get $\iota[f]=[f]$ which means that $\iota$ does not swap the summands.

Hence we have $\ker\alpha^\i_\imm\cong\Z\oplus\Z_2\oplus\Z_2$ and we can also see that $\Z$ is a direct summand in $\Imm^\i(k+1,k)$. Now if one of the $\Z_2$'s was not a direct summand in $\Imm^\i(k+1,k)$, that would mean that its generator coincided with $2[f]$ where $[f]$ is in the preimage under $\alpha^\i_\imm$ of a non-zero element in $\WW_{k+1}$, i.e. $f$ is an immersion $M^{k+1}\imto\R^{2k+1}$ where $[M]\in\WW_{k+1}$ is not zero. 
This would imply that all elements in the preimage of $[M]$ were of order $4$ or $\infty$, hence to see that both $\Z_2$'s are direct summands it is sufficient to prove that the preimage under $\alpha^\i_\imm$ of any element in $\WW_{k+1}$ 
contains an element of order $2$.

By \cite[section 4]{wallcob} we have that the $\Z_2$-algebra $\WW_*$ is generated by cobordism classes $[P(2^r-1,2^rs)]$, $[Q(2^r-1,2^rs)]$ and $[\CP^{2^r}]$ (for $r,s\ge1$) of dimensions $2^r(2s+1)-1$, $2^r(2s+1)$ and $2^{r+1}$ respectively, moreover, $P(2^r-1,2^rs)$ and $\CP^{2^r}$ are orientable while $Q(2^r-1,2^rs)$ is not cobordant to orientable manifolds. The monomials of rank $k+1$ formed by these generators are a basis of the $\Z_2$-vector space $\WW_{k+1}$. We also have $\WW_{k+1}\cong\Omega_{k+1}\oplus T_2(\Omega_k)$ and in the above basis the summands $\Omega_{k+1}$ and $T_2(\Omega_k)$ are generated respectively by the monomials containing an even and an odd number of $[Q(2^r-1,2^rs)]$'s since $\Omega_{k+1}$ contains all classes representable by orientable manifolds.

Now $\Omega_{k+1}\subset\WW_{k+1}$ is independent of the $\Z_2$'s in $\Imm^\i(k+1,k)$ which means that all basis elements of $\WW_{k+1}$ which are representable by orientable manifolds have elements of order $2$ in their $\alpha^\i_\imm$-preimage. Thus it is enough to prove that this also holds for the unorientable basis elements, that is, any monomial of rank $k+1$ formed by the manifolds $P(2^r-1,2^rs)$, $Q(2^r-1,2^rs)$ and $\CP^{2^r}$ such that it contains an odd number of $Q(2^r-1,2^rs)$'s can be immersed into $\R^{2k+1}$ representing a cobordism class of order $2$.

Observe that $Q(2^r-1,2^rs)$ is even dimensional for all $r,s\ge1$ and $k+1$ is odd, hence such a monomial can be written in the form $M^m\times N^n$ with $m,n\ge1$ and $m+n=k+1$. We can assume that $m$ is even and then \ref{iw1} and \ref{iw3} yield that there is an immersion $f\colon M^m\imto\R^{2m-1}$ such that $2[f]$ is $0$ in $\Imm^\i(m,m-1)$. But then choosing any immersion $g\colon N^n\imto\R^{2n}$ we get an immersion
$$f\times g\colon M\times N\imto\R^{2k+1}$$
which represents an element of order $2$ in $\Imm^\i(k+1,k)$ since we have $2[f\times g]=2[f]\times[g]=0$ and $\Imm^\i(*,*)$ is a bi-graded ring. This finishes our proof.
\end{prf}

\section{Morin maps}\label{morm}

Recall that a stable map $f\colon M^n\to P^{n+k}$ induces a stratification of $M$ according to its Thom--Boardman types (see \cite{boar}): $\Sigma^{i}(f)\subset M$ denotes the submanifold consisting of the points where the rank of derivative of $f$ drops by $i$, then we have $\Sigma^{i,j}(f):=\Sigma^{j}(f|_{\Sigma^i(f)})$ and so on; this defines $\Sigma^{i_1,\ldots,i_m}(f)$ for any decreasing sequence $i_1,\ldots,i_m$. In the present section we shall consider maps which only have singularities of type $\Sigma^1$, i.e. {\it Morin singularities} (besides regular germs) and we call such maps {\it Morin maps}.

Although generally the Thom--Boardman stratification is coarser than the singularity stratification, this is not the case for Morin maps, since a Morin singularity can only be of type $\Sigma^{1_r}:=\Sigma^{1,\ldots,1,0}$ (where the number of $1$'s is $r$) for some $r$ and each type $\Sigma^{1_r}$ contains precisely one singularity class; see \cite{mor}. We shall denote the singularity class in $\Sigma^{1_r}$ also by the symbol $\Sigma^{1_r}$ and note that for $r<s$ we have $\Sigma^{1_r}<\Sigma^{1_s}$, that is, Morin singularities form a single increasing sequence. We call Morin maps with at most $\Sigma^{1_r}$ singularity {\it $\Sigma^{1_r}$-maps} and put
$$\Mor_r^\sigma(n,k):=\Cob_{\{\Sigma^0,\Sigma^{1_1},\ldots,\Sigma^{1_r}\}}^\sigma(\R^{n+k})$$
for any stable normal structure $\sigma$ (except that for $\sigma=G\oplus\xi$ we decrease $n$ and increase $k$ by the rank of $\xi$). Here $r=\infty$ is also allowed and it means that we put no further restriction on Morin maps.

The rational cohomology of the Kazarian space $K_r^\SO$ of oriented $\Sigma^{1_r}$-maps (for $1\le r\le\infty$) was computed by Szűcs \cite{hosszu} and we can use this to obtain the form of $\varphi_{n+k}$ rationally:

\begin{prop}\label{morq}
The endomorphism $\varphi_{n+k}\otimes\Q$ of $\Mor_r^\SO(n,k)\otimes\Q$ is the following:
\begin{enumerate}
\item if $k=2m$ is even, $r<\infty$ is also even and $q$ is the number of non-negative integers $a_0,a_1,\ldots,a_m$ such that $a_0$ is odd and $n=a_0k(r+1)+\underset{i=1}{\overset m\sum}a_i4i$, then $\varphi_{n+k}\otimes\Q$ is trivial on $q$ generators and the multiplication by $2$ on the rest of the generators (in an appropriate basis),
\item otherwise $\varphi_{n+k}\otimes\Q$ is the multiplication by $2$.
\end{enumerate}
\end{prop}

\begin{prf}
Put $A:=\Q[p_1,\ldots,p_m]$ with $\deg p_i=4i$. By \cite[theorems 6 and 7]{hosszu} we have
\begin{alignat*}2
\Mor_r^\SO(*,k)\otimes\Q&\cong H^*(K^\SO_r;\Q)=\\
&=\begin{cases}
A,&\text{if }k=2m\text{ is even and }r\text{ is odd or }\infty\\
A[e^{r+1}]/\bigl((e^{r+1})^2-p_m^{r+1}\bigr),&\text{if }k=2m\text{ and }r\text{ are both even}\\ 
&\text{and }\deg e^{r+1}=k(r+1)\\
A,&\text{if }k=2m-1\text{ is odd and }r=\infty\\
A/\Bigl(p_m^{\lceil\frac{r+1}2\rceil}\Bigr),&\text{if }k=2m-1\text{ is odd and }r<\infty
\end{cases}
\end{alignat*}
with $p_i\in H^{4i}(K^\SO_r;\Q)$ and $e^{r+1}\in H^{k(r+1)}(K^\SO_r;\Q)$. Now following the action of the involution inducing the $\iota$ of proposition \ref{2} through the spectral sequences computed in \cite{hosszu} yields that $\iota$ changes the $p_i$ and $e^{r+1}$ here analogously to how it did in the proof of proposition \ref{immq}, that is, we get $p_i\mapsto p_i$ and $e^{r+1}\mapsto-e^{r+1}$. This implies our claim.
\end{prf}\medskip

As in the case of immersions, we only have an understanding of the torsion parts of these cobordism groups if the codimension $k$ is either small or relatively large. In particular we will investigate the cases when either we have $k=1$ and $1\le r<\infty$ or $k$ is large compared to $n$ and $r$ is $1$.

\subsection*{1-codimensional Morin maps}

Morin maps of codimension $1$ were considered by Szűcs \cite{nulladik} who computed $\Mor_r^\SO(n,1)$ for $1\le r<\infty$ modulo small torsion groups. We also note that $\Mor_r^\O(n,1)$ is finite $2$-primary for all $n$ and $1\le r\le\infty$ by \cite[theorem 1]{szszt}. In the following we denote by $\CC_2$ the Serre class of finite $2$-primary Abelian groups. 

\begin{prop}
The endomorphism $\varphi_{n+1}$ of $\Mor_r^\SO(n,1)$ is the multiplication by $2$ modulo $\CC_2$ (that is, on the odd torsion and free parts) for $1\le r<\infty$.
\end{prop}

\begin{prf}
We proceed by induction on $r$. Observe that setting $r=0$ we get the case of $1$-codimensional immersions where we saw this claim to be true in the previous section. Now assume it to be true for $r-1$ and prove for $r$.

We use the key fibration \cite[section 16]{hosszu} to obtain an exact sequence
$$\Mor_{r-1}^\SO(n,1)\to\Mor_r^\SO(n,1)\to\Imm^{\tilde\xi^\SO_r}(n-2r,2r+1)$$
where the first arrow is the natural forgetful homomorphism and the second one assigns to a cobordism class $[f]$ the cobordism class of its restriction $f|_{\Sigma^{1_r}(f)}$ to the most complicated stratum which is an immersion with normal bundle induced from $\tilde\xi^\SO_r$ (the universal normal bundle of the $\Sigma^{1_r}$-stratum of oriented maps). These both commute with the reflection to a hyperplane (without changing orientation), i.e. the involution $\iota$ in theorem \ref{t1}.

By the considerations in \cite[section 2.2]{nulladik} for $r$ odd we have $\Imm^{\tilde\xi^\SO_r}(n-2r,2r+1)\in\CC_2$ and for $r$ even the forgetful homomorphism
$$\pi^s(n-2r)\cong\Imm^\fr(n-2r,2r+1)\to\Imm^{\tilde\xi^\SO_r}(n-2r,2r+1)$$
(which assigns to the cobordism class of a framed immersion its class as an immersion with normal bundle induced from $\tilde\xi^\SO_r$) is a $\CC_2$-isomorphism. Observe that this homomorphism also commutes with the action of $\iota$.

In the previous section we saw that theorem \ref{t1} for $1$-codimensional immersions is just the long exact sequence induced by the cofibration $S^1\xra2S^1\to\RP^2$ (where $2$ denotes the degree-$2$ map), hence in that case the endomorphism $\varphi_{n-2r+1}$ of $\Imm^\SO(n-2r,1)\cong\pi^s(n-2r)\cong\Imm^\fr(n-2r,2r+1)$ is $2\id$ which means that $\iota$ acts identically on $\Imm^\fr(n-2r,2r+1)$. By the induction hypothesis $\iota$ is the identity modulo $\CC_2$ on $\Mor_{r-1}^\SO(n,1)$ too, thus it acts identically also on $\Imm^{\tilde\xi^\SO_r}(n-2r,2r+1)$ modulo $\CC_2$. But then the exact sequence above shows that $\iota$ has to be the identity modulo $\CC_2$ on $\Mor_r^\SO(n,1)$ as well which means that $\varphi_{n+1}$ is $2\id$ modulo $\CC_2$.
\end{prf}

\begin{rmk}
Using \cite[theorem A]{nulladik}, a similar argument as in the proof above also shows that for $r=1$ we have $\varphi_{n+1}=2\id$ not only modulo $\CC_2$ but even on the $2$-primary part.
\end{rmk}

\begin{crly}
For all $n$ and $1\le r<\infty$ the group $\Mor_r^\i(n,1)$ is finite $2$-primary.
\end{crly}

\begin{prf}
This is immediate from the proposition above since $\varphi_{n+1}$ is an isomorphism on the odd torsion part and the kernel and cokernel of $\varphi_{n+1}$ on the free part are $0$ and $\Z_2$ respectively.
\end{prf}\medskip

This describes $\Mor_r^\i(n,1)$ generally. Let us now consider the case $r=1$, i.e. the so-called {\it fold maps} and put $\Fold^\i(n,1):=\Mor_1^\i(n,1)$. Restricting the exact sequence of theorem \ref{t1} for fold cobordisms to the $2$-primary parts \cite[theorem A]{nulladik} yields
$$\ker(\lambda_*)_{n-1}\xra{2\id}\ker(\lambda_*)_{n-1}\to\Fold^\i(n,1)\to\ker(\lambda_*)_{n-2}\xra{2\id}\ker(\lambda_*)_{n-2}$$
where $(\lambda_*)_m$ denotes the Kahn--Priddy homomorphism $\pi^s_m(\RP^\infty)\to\pi^s(m)$ which maps onto the $2$-primary torsion part of $\pi^s(m)$ (see \cite{kp}).

Recalling the first few stable homotopy groups of $\RP^\infty$ computed by Liulevicius \cite{liu} we obtain
\begin{table}[H]\begin{center}\begin{tabular}{c||c|c|c|c|c|c|c|c|c|c}
$m$ & $0$ & $1$ & $2$ & $3$ & $4$ & $5$ & $6$ & $7$ & $8$ & $9$ \\
\hline
$\pi^s_m(\RP^\infty)$ & $0$ & $\Z_2$ & $\Z_2$ & $\Z_8$ & $\Z_2$ & $0$ & $\Z_2$ & $\Z_{16}\oplus\Z_2$ & $(\Z_2)^3$ & $(\Z_2)^4$ \\
\hline
$\pi^s(m)$ & $\Z$ & $\Z_2$ & $\Z_2$ & $\Z_8\oplus\Z_3$ & $0$ & $0$ & $\Z_2$ & $\Z_{16}\oplus\Z_3\oplus\Z_5$ & $(\Z_2)^2$ & $(\Z_2)^3$ \\
\hline
$\ker(\lambda_*)_m$ & $0$ & $0$ & $0$ & $0$ & $\Z_2$ & $0$ & $0$ & $\Z_2$ & $\Z_2$ & $\Z_2$
\end{tabular}\end{center}\end{table}\vspace{-.5cm}\noindent
using the notation $(\Z_2)^r:=\Z_2\oplus\ldots\oplus\Z_2$ with $r$ summands. This implies the following:

\begin{prop}
The cobordism groups $\Fold^\i(n,1)$ for $n\le10$ are
\begin{table}[H]\begin{center}\begin{tabular}{c||c|c|c|c|c|c|c|c|c|c|c}
$n$ & $0$ & $1$ & $2$ & $3$ & $4$ & $5$ & $6$ & $7$ & $8$ & $9$ & $10$ \\
\hline
$\Fold^\i(n,1)$ & $\Z_2$ & $0$ & $0$ & $0$ & $0$ & $\Z_2$ & $\Z_2$ & $0$ & $\Z_2$ & $\Z_2~?~\Z_2$ & $\Z_2~?~\Z_2$
\end{tabular}\end{center}\end{table}\vspace{-.5cm}\noindent
where $G~?~H$ again denotes the existence of a short exact sequence $0\to G\to G~?~H\to H\to0$.
\end{prop}

\subsection*{Large codimensional fold maps}

In the following we shall consider fold maps of codimension $k$ such that the dimension of the source manifolds is $2k+1$ or $2k+2$ and, as we also did above, use the notation
$$\Fold^\sigma(n,k):=\Mor_1^\sigma(n,k)$$
for any stable normal structure $\sigma$ (except that for $\sigma=G\oplus\xi$ we decrease $n$ and increase $k$ by the rank of $\xi$). The cobordism groups $\Fold^\SO(n,k)$ and $\Fold^\O(n,k)$ were determined in these dimensions by Ekholm, Szűcs, Terpai \cite{eszt} and Terpai \cite{2k+2}\footnote{As Terpai recently noted \cite[theorem 4.a)]{2k+2} is false. We shall elaborate more on this in remark \ref{hiba}.}.

\begin{rmk}
The analogue of remark \ref{wueiaieviaqhevahhuev} for fold maps is true for $n<2k+1$: if $n<2k+1$, then generic maps of $n$-manifolds and $(n+1)$-manifolds with codimension $k$ are stable and a computation with jet bundles implies that for any such map the codimension of the $\Sigma^i$-stratum and the $\Sigma^{1_i}$-stratum in the source manifold are $i(k+i)$ and $i(k+1)$ respectively, hence for $n<2k+1$ a generic $k$-codimensional map of a cobordism of $n$-manifolds is fold. Moreover, $n=2k+1$ and $n=2k+2$ are precisely the dimensions where fold cobordisms are not generic (that is, $\Fold^\sigma(n,k)$ is a priori not the same as the abstract cobordism group) but they are not very far from being generic as the only other type of singularity that can generically occur is the cusp, i.e. $\Sigma^{1,1,0}$.
\end{rmk}

As for the case of immersions, we will use the natural forgetful homomorphisms
\begin{gather*}
\alpha^\SO_\fold\colon\Fold^\SO(n,k)\to\Omega_n,\quad\alpha^\i_\fold\colon\Fold^\i(n,k)\to\WW_n,\\
\alpha^\O_\fold\colon\Fold^\O(n,k)\to\NN_n\quad\text{and}\quad\alpha^{\O\oplus2\Lambda\gamma}_\fold\colon\Fold^{\O\oplus2\Lambda\gamma}(n-2,k+2)\to\NN_{n-2}
\end{gather*}
that assign to the cobordism class of a fold map $f\colon M^n\to\R^{n+k}$ (or a germ $\tilde f\colon\xi_M\to\R^{n+k}$) with the appropriate normal structure the abstract cobordism class of $M$. These again commute with the exact sequences in theorems \ref{t1} and \ref{t2} and the classical exact sequences (\ref{ces1}) and (\ref{ces2}), that is, we have commutative diagrams
$$\xymatrix@C=.75pc{
\ldots\ar[r] & \Fold^\SO(n,k)\ar[d]^{\alpha^\SO_\fold}\ar[r] & \Fold^\SO(n,k)\ar[d]^{\alpha^\SO_\fold}\ar[r] & \Fold^\i(n,k)\ar[d]^{\alpha^\i_\fold}\ar[r] & \Fold^\SO(n-1,k)\ar[d]^{\alpha^\SO_\fold}\ar[r] & \ldots \\
\ldots\ar[r] & \Omega_n\ar[r]^{2\id} & \Omega_n\ar[r] & \WW_n\ar[r] & \Omega_{n-1}\ar[r]^{2\id} & \ldots
}$$
and
$$\xymatrix@C=.75pc{
\ldots\ar[r] & \Fold^\i(n,k)\ar[d]^{\alpha^\i_\fold}\ar[r] & \Fold^\O(n,k)\ar[d]^{\alpha^\O_\fold}\ar[r] & \Fold^{\O\oplus2\Lambda\gamma}(n-2,k+2)\ar[d]^{\alpha^{\O\oplus2\Lambda\gamma}_\fold}\ar[r] & \Fold^\i(n-1,k)\ar[d]^{\alpha^\i_\fold}\ar[r] & \ldots \\
\ldots\ar[r]^0 & \WW_n\ar[r] & \NN_n\ar[r] & \NN_{n-2}\ar[r]^0 & \WW_{n-1}\ar[r] & \ldots
}$$
where the rows are exact.

\begin{lemma}\label{foldo+}
The map $\alpha^{\O\oplus2\Lambda\gamma}_\fold$ is an isomorphism $\Fold^{\O\oplus2\Lambda\gamma}(n-2,k+2)\cong\NN_{n-2}$ for all $n\le2k+2$.
\end{lemma}

\begin{prf}
This is almost completely analogous to lemma \ref{immo+} by changing immersions to fold maps and $k+1$ to $2k+2$. The only thing we have to add to its proof now is that the generic map
$$F\colon\xi_W\to\R^{n+k}\times[0,1]$$
(where $W$ is a cobordism of $(n-2)$-manifolds and $\xi_W:=2\Lambda TW$) only has fold singularities in a neighbourhood of the zero-section $W$ and its germ along $W$ satisfies the condition in definition \ref{deftauxi} that the differential $d\tilde F$ restricted to any fibre of $\xi_W$ is injective.

We may assume that the singularity strata intersect $W\subset\xi_W$ transversally which firstly means that $W$ is disjoint from the cusp stratum (since the cusp stratum is at most $1$-dimensional and $W$ is $2$-codimensional in $\xi_W$), hence it is indeed a fold map on a small neighbourhood of $W$. Secondly it means that the local trivialisations of $\xi_W$ in a sufficiently small neighbourhood of the zero-section can be chosen such that for all $p\in\Sigma^{1,0}(F)$ the whole fibre (in this neighbourhood) of $\xi_W$ containing $p$ belongs to $\Sigma^{1,0}(F)$, that is, both the fibres of $\xi_W$ and the fold stratum are orthogonal to $W$. 

Now for any point $p\in\Sigma^{1,0}(F)\cap W$ a coordinate neighbourhood of $p$ has the form $\R^{k+1}\times\R^{n-k-2}\times\R^2$ where $\R^{k+1}\times\R^{n-k-2}$ and $\R^{n-k-2}\times\R^2$ are coordinate neighbourhoods in $W$ and in $\Sigma^{1,0}(F)$ respectively. By \cite[theorem 6]{rsz} the normal bundle of $\Sigma^{1,0}(F)$ is induced from $\gamma_1^\O\oplus\gamma_1^\O\otimes\gamma_k^\O$ over $BG_1^\SO$ with $G_1^\SO=\{(\varepsilon,A)\in\O(1)\times\O(k)\mid\varepsilon\det A=1\}$. The fibres of this in the coordinate neighbourhood of $p$ can be assumed to be tangent to $\R^{k+1}\times\{x\}$ for all $x\in\R^{n-k-2}\times\R^2$ and can also be identified over all points $x$. Thus in this neighbourhood the map $F$ has the desired form described in definition \ref{deftauxi}. If $p$ was instead in $W\setminus\Sigma^{1,0}(F)$ then $F$ is an immersion near $p$, hence then we also get the desired form of $F$ and so the first condition in definition \ref{deftauxi} is always satisfied.
\end{prf}

\begin{prop}\label{foldi}
For all $k\ge1$ we have
\begin{enumerate}
\item\label{w1} $\Fold^\i(2k+1,k)\cong\WW_{2k+1}$ if $k\ne2$,
\item\label{w2} $\Fold^\i(5,2)\cong\WW_5\oplus\Z_2\cong\Z_2\oplus\Z_2$.
\end{enumerate}
\end{prop}

\begin{prf}
We shall use the first diagram above lemma \ref{foldo+} and obtain $\Fold^\SO(2k+1,k)$ from \cite[theorem 1]{eszt}.

\medskip\noindent\textit{Proof of \ref{w1}.\enspace\ignorespaces}
If $k>2$ is even, then $\alpha^\SO_\fold$ is an isomorphism in this dimension yielding the diagram
$$\xymatrix{
\Omega_{2k+1}\ar[d]^{\id}\ar[r]^{2\id} & \Omega_{2k+1}\ar[d]^{\id}\ar[r] & \Fold^\i(2k+1,k)\ar[d]^{\alpha^\i_\fold}\ar[r] & \Omega_{2k}\ar[d]^{\id}\ar[r]^{2\id} & \Omega_{2k}\ar[d]^{\id} \\
\Omega_{2k+1}\ar[r]^{2\id} & \Omega_{2k+1}\ar[r] & \WW_{2k+1}\ar[r] & \Omega_{2k}\ar[r]^{2\id} & \Omega_{2k}
}$$
and now the $5$-lemma implies the isomorphism claimed.

For $k$ odd we have $\Fold^\SO(2k+1,k)\cong\Omega_{2k+1}\oplus\Z_{3^t}$ (with an appropriate number $t$) and this gives
$$\xymatrix{
\Omega_{2k+1}\oplus\Z_{3^t}\ar[r]^{\id+\iota} & \Omega_{2k+1}\oplus\Z_{3^t}\ar[d]^{\id\oplus0}\ar[r] & \Fold^\i(2k+1,k)\ar[d]^{\alpha^\i_\fold}\ar[r] & \Omega_{2k}\ar[d]^{\id}\ar[r] & 0 \\
0\ar[r] & \Omega_{2k+1}\ar[r] & \WW_{2k+1}\ar[r] & \Omega_{2k}\ar[r] & 0
}$$
thus $\ker\alpha^\i_\fold$ is $\Z_{2^u}$ for some $0\le u\le t$ depending on how $\iota$ acts on the $\Z_{3^t}$ part of $\Fold^\SO(2k+1,k)$. Now applying lemma \ref{foldo+}, part of the second diagram above it has the form
$$\xymatrix{
\NN_{2k}\ar[r]^(.33){\psi'_{3k+2}} & \Fold^\i(2k+1,k)\ar[d]^{\alpha^\i_\fold}\ar[r]^(.67){\varphi'_{3k+1}} & \NN_{2k+1}\ar[d]^{\id} \\
0\ar[r] & \WW_{2k+1}\ar[r] & \NN_{2k+1}
}$$
hence we have $\Z_{3^u}\cong\ker\alpha^\i_\fold=\ker\varphi'_{3k+1}=\im\psi'_{3k+2}$ which is isomorphic to a factor group of $\NN_{2k}$. But $\NN_{2k}$ is $2$-primary so this can only happen if $u=0$, that is, if $\alpha^\i_\fold$ is an isomorphism.

\medskip\noindent\textit{Proof of \ref{w2}.\enspace\ignorespaces}
We have $\Fold^\SO(4,2)\cong\Omega_4\cong\Z$ and $\Fold^\SO(5,2)\cong\Omega_5\oplus\Z_2\cong\Z_2\oplus\Z_2$ and their endomorphisms $\varphi_6$ and $\varphi_7$ are both the multiplication by $2$, hence we get the diagram
$$\xymatrix{
0\ar[r] & \Omega_5\oplus\Z_2\ar[d]^{\id\oplus0}\ar[r] & \Fold^\i(5,2)\ar[d]^{\alpha^\i_\fold}\ar[r] & 0 \\
0\ar[r] & \Omega_5\ar[r] & \WW_5\ar[r] & 0
}$$
which proves our statement.
\end{prf}

\begin{prop}
For all $k\ge1$ we have
\begin{enumerate}
\item\label{fw1} $\Fold^\i(2k+2,k)$ is isomorphic to a subgroup of $\WW_{2k+2}$ of index $2$ if $k\ne2$,
\item\label{fw2} $\Fold^\i(6,2)\cong\WW_6\oplus\Z_2\cong\Z_2\oplus\Z_2$.
\end{enumerate}
\end{prop}

\begin{prf}
We will again always use the first diagram above lemma \ref{foldo+} and we know $\Fold^\SO(2k+1,k)$ and $\Fold^\SO(2k+2,k)$ from \cite[theorem 1]{eszt} and \cite[theorem 4]{2k+2} respectively. 

\medskip\noindent\textit{Proof of \ref{fw1}.\enspace\ignorespaces}
If $k>2$ is even, then $\alpha^\SO_\fold$ is an isomorphism $\Fold^\SO(2k+1,k)\cong\Omega_{2k+1}$ and the embedding of $\Fold^\SO(2k+2,k)$ into $\Omega_{2k+2}$ as an index-$2$ subgroup.

If $k$ is odd, then we have $\Fold^\SO(2k+1,k)\cong\Omega_{2k+1}\oplus\Z_{3^t}$ but, as we saw in the proof of proposition \ref{foldi}, the $\Z_{3^t}$ part is mapped isomorphically by the endomorphism $\varphi_{3k+1}$, hence it does not appear in its kernel. Moreover, for $k$ odd $\alpha^\SO_\fold$ embeds $\Fold^\SO(2k+2,k)$ into $\Omega_{2k+2}$ as the kernel of the homomorphism $\overline p_{\frac{k+1}2}[\cdot]\colon\Omega_{2k+2}\to\Z$ which maps a cobordism class $[M]$ to the normal Pontryagin number $\overline p_{\frac{k+1}2}[M]$. Now if we take the cokernel of $\varphi_{3k+2}$, that is, we factor $\ker\overline p_{\frac{k+1}2}[\cdot]$ by $2\ker\overline p_{\frac{k+1}2}[\cdot]$ and $\Omega_{2k+2}$ by $2\Omega_{2k+2}$, then the image of $\Fold^\SO(2k+2,k)$ under this quotient map will again become an index-$2$ subgroup.

Hence both for $k$ odd and for $k$ even we obtain that there is a group $G$ such that the commutative diagram
$$\xymatrix{
0\ar[r] & G\ar@{^(->}[d]\ar[r] & \Fold^\i(2k+2,k)\ar[d]^{\alpha^\i_\fold}\ar[r] & \Omega_{2k+1}\ar[d]^\id\ar[r] & 0 \\
0\ar[r] & \frac{\Omega_{2k+2}}{2\Omega_{2k+2}}\ar[r] & \WW_{2k+2}\ar[r] & \Omega_{2k+1}\ar[r] & 0
}$$
holds where the monomorphism from $G$ into $\Omega_{2k+2}/2\Omega_{2k+2}$ has cokernel $\Z_2$. This implies our statement.

\medskip\noindent\textit{Proof of \ref{fw2}.\enspace\ignorespaces}
We have $\Fold^\SO(5,2)\cong\Omega_5\oplus\Z_2\cong\Z_2\oplus\Z_2$ and $\Fold^\SO(6,2)\cong\Omega_6\cong0$ which gives us the diagram
$$\xymatrix{
0\ar[r] & \Fold^\i(6,2)\ar[d]^{\alpha^\i_\fold}\ar[r] & \Omega_5\oplus\Z_2\ar[d]^{\id\oplus0}\ar[r] & 0 \\
0\ar[r] & \WW_6\ar[r] & \Omega_5\ar[r] & 0
}$$
proving our claim.
\end{prf}


\begin{crly}
For all $k\ne2$ the homomorphism $\alpha^\O_\fold$ embeds $\Fold^\O(2k+2,k)$ into $\NN_{2k+2}$ as a subgroup of index $2$.
\end{crly}

\begin{prf}
We apply lemma \ref{foldo+}, the second diagram above it and also that the homomorphism $\psi'_{3k+3}\colon\Fold^{\O\oplus2\Lambda\gamma}(2k+1,k+2)\to\Fold^\i(2k+2,k)$ is zero as we saw in the proof of proposition \ref{foldi}. Then we have
$$\xymatrix{
& \Fold^\i(2k+2,k)\ar[d]^{\alpha^\i_\fold}\ar[r] & \Fold^\O(2k+2,k)\ar[d]^{\alpha^\O_\fold}\ar[r] & \NN_{2k}\ar[d]^\id\ar[r] & 0 \\
0\ar[r] & \WW_{2k+2}\ar[r] & \NN_{2k+2}\ar[r] & \NN_{2k}\ar[r] & 0
}$$
thus the proposition above implies that $\alpha^\O_\fold$ is a monomorphism with cokernel $\Z_2$.
\end{prf}

\begin{rmk}\label{hiba}
This corollary contradicts \cite[theorem 4.a)]{2k+2} which essentially claims that both the kernel and the cokernel of
$$\alpha^\O_\fold\colon\Fold^\O(2k+2,k)\to\NN_{2k+2}$$
are $\Z_2$. The problem with this, according to Terpai, is the following: the computations in \cite{2k+2} show that here $\ker\alpha^\O_\fold$ appears as the factor by $\Z_2$ of the group $\Z_2\oplus\Z_2$ of a complete set of invariants of cusp null-cobordisms of fold maps (recall that a null-cobordism is a map of a $(2k+3)$-manifold with codimension $k$ and such a map generically only has $\Sigma^{1,0}$ (fold) and $\Sigma^{1,1,0}$ (cusp) singularities); one of the summands $\Z_2$ (which gets cancelled by the factoring) measures the number of components of the cusp curve of such a null-cobordism modulo $2$ and the other one measures whether the orientation changes over the cusp curve or not. However, if we apply \cite[lemma 2]{2k+2} correctly, it yields that this second $\Z_2$ does not appear as the orientation cannot change over this curve. Hence $\ker\alpha^\O_\fold$ is actually the factor of only $\Z_2$ by $\Z_2$, that is, $\alpha^\O_\fold$ is mono. Moreover, it is also apparent that its image is the kernel of the Thom polynomial of cusp singularity $(\ol w_{k+1}^2+\ol w_k\ol w_{k+2})[\cdot]\colon\NN_{2k+2}\to\Z_2$ which maps a cobordism class $[M]$ to the normal Stiefel--Whitney number $\ol w_{k+1}^2[M]+\ol w_k\ol w_{k+2}[M]$.
\end{rmk}




\section{Bordisms}\label{bor}

In this final section we point out that theorems \ref{t1} and \ref{t2} also hold when we change $\tau$-cobordism groups to general bordism groups (although the set of all possible singularities of $k$-codimensional germs generally contains unstable ones which we excluded from $\tau$). This can most easily be seen by using somewhat more refined forms of the same arguments as in remarks \ref{gen1} and \ref{gen2}. 
This way we get the following:

\begin{thm}\label{b}
For any CW complex $Q$ there is a long exact sequence
$$\ldots\to\Omega_n(Q)\xra{\varphi_n}\Omega_n(Q)\xra{\chi_n}\WW_n(Q)\xra{\psi_n}\Omega_{n-1}(Q)\to\ldots$$
and a short exact sequence
$$0\to\WW_n(Q)\xra{\varphi'_n}\NN_n(Q)\xra{\chi'_n}\NN_{n-2}(Q)\to0$$
where $\varphi_n$ is the multiplication by $2$, $\chi_n$ and $\varphi'_n$ are the forgetful homomorphisms, $\psi_m$ assigns to the bordism class of the map of a Wall manifold $f\colon M\to Q$ the bordism class of $f|_{\PD(w_1(M))}$ and $\chi'_n$ assigns to the bordism class of a map $g\colon N\to Q$ the bordism class of $g|_{\PD(w_1(N)^2)}$. 
\end{thm}

\begin{prf}
The $n$'th bordism groups of any CW complex are the same as those of its $(n+1)$-skeleton, hence for the purpose of this proof we may only take a finite dimensional skeleton of $Q$ instead of the actual $Q$. A finite dimensional CW complex is homotopy equivalent to an orientable manifold (by embedding it into a Eucledian space and taking a small neighbourhood of it there), and so we may also assume that $Q$ is an orientable manifold of dimension $q$.

Fix a number $n$ and let $r$ be such that generic maps of $(n+1)$-manifolds to $r$-manifolds are embeddings. Then we have
$$\Omega_n(Q)\cong\Cob_\tau^\SO(Q\times\R^r),\quad\WW_n(Q)\cong\Cob_\tau^\i(Q\times\R^r)\quad\text{and}\quad\NN_n(Q)\cong\Cob_\tau^\O(Q\times\R^r)$$
for any set $\tau$ of stable singularities of codimension $q+r-n$. 

\medskip\begin{sclaim}
$\NN_{n-2}(Q)\cong\Cob_\tau^{\O\oplus2\Lambda\gamma}(Q\times\R^r)$.
\end{sclaim}

\begin{sprf}
Let $f\colon M^{n-2}\to Q$ be any map of an $(n-2)$-manifold to $Q$ and put $\xi_M:=2\Lambda\nu_f$
. Now let $\tilde f\colon\xi_M\to Q$ be any map that extends the map $f$ given on the zero-section of $\xi_M$ (such an extension is unique up to homotopies fixed on the zero-section) and let $i\colon\xi_M\into\R^r$ be any embedding (which is also unique up to isotopy). Then assigning to the map $f\colon M\to Q$ the germ of the embedding
$$\tilde f\times i\colon\xi_M\into Q\times\R^r$$
along the zero-section is well-defined on bordism classes up to cobordism of embeddings with normal $\O\oplus2\Lambda\gamma$-structure, hence it yields a homomorphism
$$\alpha\colon\NN_{n-2}(Q)\to\Cob_\tau^{\O\oplus2\Lambda\gamma}(Q\times\R^r).$$

The inverse of this $\alpha$ can be constructed as follows: a map germ $\xi_M\into Q\times\R^r$ is generically the germ of an embedding of the form $\tilde f\times i$ where $\tilde f$ is any map to $Q$ and $i$ is an embedding into $\R^r$, and then we can assign to the cobordism class of the germ of $\tilde f\times i$ the bordism class of the restriction $\tilde f|_M$. This is again a well-defined homomorphism and it is not hard to see that it inverts $\alpha$, thus $\alpha$ is an isomorphism. 
\end{sprf}

The above considerations and theorems \ref{t1} and \ref{t2} give two exact sequences
$$\ldots\to\Omega_n(Q)\xra{\varphi_n}\Omega_n(Q)\xra{\chi_n}\WW_n(Q)\xra{\psi_n}\Omega_{n-1}(Q)\to\ldots$$
and
$$\ldots\to\WW_n(Q)\xra{\varphi'_n}\NN_n(Q)\xra{\chi'_n}\NN_{n-2}(Q)\xra{\psi'_n}\WW_{n-1}(Q)\to\ldots$$
infinite to the right, moreover, keeping $n$ fixed and growing $r$ yields that they are also infinite to the left. By design the homomorphisms $\chi_n$, $\psi_n$, $\varphi'_n$ and $\chi'_n$ are as in theorems \ref{t1} and \ref{t2}. 

\medskip\begin{sclaim}
The endomorphism $\varphi_n$ of $\Omega_n(Q)$ is the multiplication by $2$.
\end{sclaim}

\begin{sprf}
We want to see that the involution $\iota$ of $\Omega_n(Q)\cong\Cob_\tau^\SO(Q\times\R^r)$ is the identity, that is, composing an oriented embedding $e\colon M^n\into Q\times\R^r$ with the reflection to a hyperplane but keeping its orientation unchanged will be cobordant to $e$. Such a reflection does not change the homotopy class of $e$, thus there is a map
$$E\colon M\times[0,1]\to Q\times\R^r\times[0,1]$$
(mapping $M\times\{t\}$ to $Q\times\R^r\times\{t\}$ for all $t\in[0,1]$) extending $e$ on $M\times\{0\}$ and its reflection on $M\times\{1\}$. Moreover, $\nu_E$ is canonically oriented since the orientations given on the normal bundles of $e$ and its reflection are identical. Now since $E$ is generically an embedding, it is an oriented cobordism between $e$ and its reflection meaning that we have $\iota[e]=[e]$. Thus we got $\varphi_n=2\id$ as we wanted.
\end{sprf}

\begin{sclaim}
The homomorphism $\psi'_n\colon\NN_{n-2}(Q)\to\WW_{n-1}(Q)$ is zero.
\end{sclaim}

\begin{sprf}
It is enough to show a splitting $\sigma_n\colon\NN_{n-2}(Q)\to\NN_n(Q)$ of the sequence above. This is exactly the same as in \cite[theorem 4.4]{atiyah}, namely we define $\sigma_n$ to assign to a bordism class represented by a map $f\colon M^{n-2}\to Q$ the bordism class of $f\circ\pi$ where $\pi$ is the projection of the $\RP^2$-bundle $P(\Lambda TM\oplus\varepsilon^2)$ (the projectivisation of $\Lambda TM\oplus\varepsilon^2$) over $M$. The proof in Atiyah's paper also works in this setting and gives that $\sigma_n$ is well-defined and $\chi'_n\circ\sigma_n$ is the identity.
\end{sprf}


Now the fact that $\varphi'_n\colon\WW_n(Q)\to\NN_n(Q)$ is mono implies that in the definition of $\psi_n\colon\WW_n(Q)\to\Omega_{n-1}(Q)$ we do not have to differentiate between the integer representatives of the $w_1$ class. Thus we have proved everything we wanted.
\end{prf}

\begin{rmk}
~\vspace{-.5em}

\begin{enumerate}
\item The first exact sequence in theorem \ref{b} was also obtained by Stong \cite[proposition 6.1]{stong} and partly by Conner and Floyd \cite[1.10 theorem]{cf}.
\item If the $Q$ in theorem \ref{b} is an (unorientable) manifold, then the theorem also holds when we change $\Omega_n(Q)$, $\WW_n(Q)$, $\psi_n$ and $\chi'_n$ to their normal analogues, i.e. changing everywhere the tangent $w_1$ class to the normal $w_1$ class.
\end{enumerate}
\end{rmk}

\begin{rmk}
Part of the results collected in theorems \ref{t1}, \ref{t2} and \ref{b} can also be proved using methods analogous to those of Wall \cite{wallcob} and Li \cite{li}; the reader is encouraged to try working them out for themselves.
\end{rmk}

Lastly we note that the second exact sequence in theorem \ref{b} can be obtained directly from the classical exact sequence (\ref{ces2}) because of the following (we note here that this proposition is essentially due to Terpai):

\begin{prop}\label{ww=}
For any CW complex $Q$ we have
$$\WW_n(Q)=\bigoplus_{i+j=n}\WW_i\otimes H_j(Q;\Z_2).$$
\end{prop}

\begin{prf}
We imitate the proof of the fact
$$\NN_n(Q)=\bigoplus_{i+j=n}\NN_i\otimes H_j(Q;\Z_2);$$
see \cite{cf}. Let $E_{i,j}^m$ be the Atiyah--Hirzebruch spectral sequence for Wall bordisms, i.e. whose $E^\infty$-page is associated to a filtration of $\WW_*(Q)$. Then we have $E^2_{i,j}=H_i(Q;\Z_2)\otimes\WW_j$ since $\WW_n$ is a $\Z_2$-vector space.

Now it is sufficient to see that the spectral sequence is trivial. By \cite{cf} this is equivalent to saying that the homomorphism $\mu_n\colon\WW_n(Q)\to H_n(Q;\Z_2)$ which maps a bordism class to the homology class represented by its elements is epi. In other words we want to show that any $\Z_2$-homology class of $Q$ is representable by the map of a Wall manifold to $Q$.

Thom \cite{th} showed that any homology class $x\in H_n(Q;\Z_2)$ can be represented by a map $f\colon M^n\to Q$ where $M$ is a manifold; our goal is to show that $M$ can also be chosen to have integer $w_1$ class. Choose a submanifold $N^{n-2}\subset M^n$ representing $\PD(w_1(M)^2)$ and obtain $b_N\colon B(N)\to M$ by blowing up $N$, that is, for all points $p\in N$ we change the normal disk of $N$ at $p$ to a Möbius strip; $b_N^{-1}(p)=S^1$ is the central curve of this Möbius strip and otherwise $b_N|_{B(N)\setminus b_N^{-1}(N)}$ is a diffeomorphism to $M\setminus N$. Then $f\circ b_N$ also represents the homology class $x$. We claim that $B(N)$ has integer first Stiefel--Whitney class which implies our statement.

A more precise description of $B(N)$ is the following: the normal bundle of $N\subset M$ is $2\delta$ where $\delta$ is the determinant line bundle $\det TN=\det TM|_N$; now $N$ is also a submanifold of the $\RP^2$-bundle $P(\delta\oplus\varepsilon^2)$ with normal bundle $2\delta$, hence the disk bundle $D(2\delta)$ appears as a tubular neighbourhood of $N$ in both $M$ and $P(\delta\oplus\varepsilon^2)$. We have
$$B(N)=(M\setminus D(2\delta))\usqcup{S(2\delta)}(P(\delta\oplus\varepsilon^2)\setminus D(2\delta))$$
where this denotes the gluing along the boundaries $S(2\delta)$.

To show that $w_1(B(N))$ is integer 
it is enough to show that the Poincaré dual of its square can be represented by the empty submanifold. The reason for this is that on one hand the integer $w_1$ class means that $w_1(B(N))\colon B(N)\to\RP^\infty$ maps to the $1$-cell and on the other hand $\PD(w_1(B(N))^2)=\varnothing$ means that the image of $w_1(B(N))$ is disjoint from $\RP^{\infty-2}$ and these two conditions are equivalent up to homotopy.

We have $\PD(w_1(B(N))^2)=\PD(w_1(M\setminus D(2\delta))^2)\cup\PD(w_1(P(\delta\oplus\varepsilon^2)\setminus D(2\delta))^2)$ when represented by manifolds matching on the boundaries. Now $\PD(w_1(M)^2)$ is represented by the intersection $N$ of two transverse submanifolds $V$ and $V'$ representing $\PD(w_1(M))$; then $\PD(w_1(M\setminus D(2\delta)))=\PD(w_1(M)|_{M\setminus D(2\delta)})$ is represented by $V\cap(M\setminus D(2\delta))=V\setminus D(2\delta)$ and by $V'\cap(M\setminus D(2\delta))=V'\setminus D(2\delta)$ which means that their intersection representing $\PD(w_1(M\setminus D(2\delta))^2)$ is empty. But we know from the proof of \cite[theorem 4.4]{atiyah} that $\PD(w_1(P(\delta\oplus\varepsilon^2))^2)$ can be represented by $N$ embedded as the zero section of $D(2\delta)$, hence its intersection with the complement of $D(2\delta)$ representing $\PD(w_1(P(\delta\oplus\varepsilon^2)\setminus D(2\delta))^2)$ is also empty. This finishes the proof.
\end{prf}\medskip

 The method of the above proof implies that the forgetful homomorphism $\varphi'_n\colon\WW_n(Q)\to\NN_n(Q)$ respects the decompositions
$$\WW_n(Q)=\bigoplus_{i+j=n}\WW_i\otimes H_j(Q;\Z_2)\quad\text{and}\quad\NN_n(Q)=\bigoplus_{i+j=n}\NN_i\otimes H_j(Q;\Z_2).$$
Thus the second exact sequence in theorem \ref{b} is the degree-$n$ part of the tensor product of the sequence (\ref{ces2})
$$0\to\WW_*\to\NN_*\to\NN_{*-2}\to0$$
with $H_*(Q;\Z_2)$.













\medskip

\begin{ack}
I thank A. Szűcs and T. Terpai for the several consultations we had on the topic of this paper. They both made essential contributions: the main ideas in sections \ref{tauw} and \ref{lessect} are due to Szűcs; Terpai resolved a contradiction in section \ref{morm}; and section \ref{bor} is mostly based on one consultation we had with Szűcs and Terpai.
\end{ack}


\end{document}